\documentclass[11pt,a4paper]{article}
\usepackage[bottom=2.0cm,top=2.0cm,left=2.0cm,right=2.0cm]{geometry}
\usepackage{graphicx} 
\usepackage{mathrsfs}
\usepackage[english]{babel}
\usepackage{calrsfs}
\usepackage{subcaption}
\usepackage{authblk}
 
\usepackage{verbatim}
\usepackage{tikz}
\usepackage{footmisc}
\usetikzlibrary{positioning}
\usepackage{pgfplots}
\usetikzlibrary{external}
\usetikzlibrary{shapes}
\usetikzlibrary{positioning,arrows.meta}
\pgfplotsset{compat=newest} 
\usepgfplotslibrary{units} 
\usepackage{amsfonts}
\usepackage{booktabs}
\usepackage{hyperref}
\usepackage{amssymb}
\usepackage{amsmath}
\usepackage{amsthm}

\usepackage{subcaption}
\usepackage{xcolor}
\usepackage{cases}
\usepackage{bm}
\usepackage[overload]{empheq} 
\usepackage{fix-cm} 

\usepackage{empheq,etoolbox}
\usepackage{float}
\usepackage{matlab-prettifier}
\usepackage{epstopdf}
\usepackage{inputenc}

\usepackage{textcomp}
\usepackage{tabularx}
\usepackage{longtable} 
\usepackage{colortbl}
\usepackage[justification=centering]{caption}

\usepackage{algorithm}
\usepackage{algorithmic}

\usepackage{appendix}

\def\tsc#1{\csdef{#1}{\textsc{\lowercase{#1}}\xspace}}
\tsc{WGM}
\tsc{QE}

\newtheorem{rmk}{Remark}
\newtheorem{assumption}{Assumption}

\newcommand{\averagel}{\{\!\!\{}
\newcommand{\averager}{\}\!\!\}}

\newcommand{\jumpl}{[\![}
\newcommand{\jumpr}{]\!]}

\newcommand{\jjumpl}{[\![\![}
\newcommand{\jjumpr}{]\!]\!]}
\DeclareMathOperator{\tr}{tr}

\mathtoolsset{showonlyrefs}
\usepackage{ulem}

\mathtoolsset{showonlyrefs,showmanualtags}

\title{A coupled mathematical and numerical model for protein spreading and tissue atrophy, applied to Alzheimer’s disease \footnote{\textbf{Funding}: VP, MC, and PFA have been funded by the European Union (ERC, NEMESIS, project number 101115663). Views and opinions expressed are however those of the author(s) only and do not necessarily reflect those of the European Union or the European Research Council Executive Agency. The present research is part of the activities of Dipartimento di Eccellenza 2023-2027. DR has been partially supported by PRIN 2022 project \textit{Mathematical models for viscoelastic biological matter}, Prot. 202249PF73 – Funded by European Union - Next Generation EU - Italian Recovery and Resilience Plan (PNRR) - M4C1\_CUP D53D23005610001. MC and PFA are members of INdAM-GNCS.DR is member of INdAM-GNFM.}}

\author[1]{Valentina Pederzoli\footnote{valentina.pederzoli@polimi.it}}

\author[1]{Mattia Corti\footnote{mattia.corti@polimi.it}}

\author[1]{Davide Riccobelli\footnote{davide.riccobelli@polimi.it}}

\author[1]{Paola F. Antonietti\footnote{paola.antonietti@polimi.it}}

\affil[1]{MOX-Dipartimento di Matematica, Politecnico di Milano, Piazza Leonardo da Vinci 32, Milan, 20133, Italy}

\date{December 2024}

\begin{document}
\maketitle
\begin{abstract}
    The aim of this paper is to introduce, analyse and test in practice a new mathematical model describing the interplay between biological tissue atrophy driven by pathogen diffusion, with applications to neurodegenerative disorders. 
    This study introduces a novel mathematical and computational model comprising a Fisher-Kolmogorov equation for species diffusion coupled with an elasticity equation governing mass loss. These equations intertwine through a logistic law dictating the reduction of the medium's mass.
    One potential application of this model lies in understanding the onset and development of Alzheimer's disease. Here, the equations can describe the propagation of misfolded $\tau$-proteins and the ensuing brain atrophy characteristic of the disease.
    To address numerically the inherited complexities, we propose a Polygonal Discontinuous Galerkin method on polygonal/polyhedral grids for spatial discretization, while time integration relies on the $\theta$-method. We present the mathematical model, delving into its characteristics and propose discretization applied.
    Furthermore, convergence results are presented to validate the model, accompanied by simulations illustrating the application scenario of the onset of Alzheimer's disease.
\end{abstract}
\section{Introduction}
\label{sec:introduction}

Population dynamics is increasingly employed in biology to describe the evolution of various phenomena, including the progression of neurodegenerative diseases \cite{physicsbasedmodelAlz} and the development of tumors \cite{FKtumor}. In some cases, the spread of pathogens can cause changes in mass within the affected tissues. The integration of population dynamics systems with models of growth or atrophy allows the description of such phenomena employing coupled multiphysics models. For instance, neurodegenerative diseases arise from the damage and degeneration of the neurons in the regions of the brain associated with cognitive functions \cite{juckerSelfpropagationPathogenicProtein2013}. 
\par
In the context of neurodegenerative diseases the underlying causes often involve the gradual buildup of damaged protein agglomerations, which cause the eventual degeneration of neurons. After sustaining damage, it is theorized that the protein gains the capacity to migrate between neurons, induces misfolding of other healthy proteins, and aggregates with them \cite{ActaNeuro,propspread, selfprop}. This process is known as prion-like behavior and is used to model the spreading of the misfolded proteins. In the literature, we can find different approaches to model the spreading of misfolded proteins in the framework of neurodegenerative diseases. One major class includes: kinetic growth and fragmentation models, which employ a set of ordinary differential equations to study the local interaction of aggregates of different sizes \cite{schiesser_parkinson}. Alternative approaches employ, network diffusion models and graph theory, used to study the global prion-like spreading of misfolded proteins \cite{Prionlike}. A third widely employed class of approaches builds on reaction-diffusion continuum models, employing systems of partial differential equations, to study the spatiotemporal evolution of the concentration of the misfolded proteins \cite{physicsbasedmodelAlz}. In the literature, three possible models are proposed: the Fisher-Kolmogorov equation, the Heterodimer model, and the Smoluchowski model \cite{Prionlike, physicsbasedmodelAlz}. The Fisher-Kolmogorov model is a nonlinear reaction-diffusion equation in one variable modeling, the relative concentration of the pathogen \cite{physicsbasedmodelAlz, Prionlike, FK_equation_and_PolyDG}. It is widely used for its simplicity even though it can not capture the mechanism of infection or the intermediate states. The heterodimer model is more complex as it accounts for two different configurations of proteins: healthy and misfolded ones. In this model, the rates of aggregation and conversion of the proteins are taken into account, as well as the clearance and production of misfolded and healthy proteins \cite{physicsbasedmodelAlz, Heterodimer}. However, this model cannot capture the size of the misfolded proteins aggregate, or their nucleation and fragmentation \cite{Prionlike}. Finally, the Smoluchowski model is the more complex approach but allows the study the kinetics of protein aggregates of different sizes \cite{smoluchowski_kuhl}, while involving a large number of parameters to be calibrated \cite{Prionlike}. To study how the spreading of misfolded proteins influences the atrophy of the brain during the development of neurodegenerative diseases, these propagation models can be coupled with a model of atrophy, where the loss of mass is affected by the concentration of misfolded proteins. The coupling is achieved by introducing a measure of the volume loss, the relative rate of which is somewhat proportional to the total exposure of the tissue to the misfolded proteins \cite{biochemical_biomechanical_Kuhl,multiph}.
\par
We introduced a multiphysics model to investigate tissue atrophy driven by pathogen spread, with a focus on applications to Alzheimer's disease. The novelty of the model relies on a new constitutive equation for the inelastic component of the stress tensor used to model tissue atrophy. This method, based on the usage of polytopal elements, facilitates the discretization of complex domains.

The propagation and aggregation of the pathogen are modeled using the Fisher-Kolmogorov equation, capturing the effects of both dispersion and proliferation. The FK equation outlines the evolution of population density, i.e., the pathogen concentration under analysis. It serves as the simplest reaction-diffusion equation that incorporates two critical effects: dispersion (diffusive term) and proliferation (reactive term) \cite{FKtumor}. Tissue atrophy is characterized through a morpho-elastic framework, which combines the effects of mass loss and tissue elasticity to determine the resulting tissue morphology, employing a multiplicative decomposition of the deformation gradient into an elastic and a growth-related component \cite{Goriely_Math_Mec}, where the latter depends on the evolution of tissue loss \cite{Ambrosimorpho}. The model connects the morpho-elastic response to pathogen concentration by defining an evolution law for inelastic strain, regulated by pathogen concentration through a logistic-type differential equation.
For the numerical discretization of the resulting coupled problem, we propose to employ a high-order PolyDG formulation for the spatial discretization coupled with the $\theta$ method for time integration, and we consider a semi-implicit approach to treat the nonlinear terms in the FK equation and in the coupling between species concentration and medium mass reduction.
\par
The numerical implementations found in literature commonly use a continuous finite element method (FEM). Another approach can be found in \cite{FK_equation_and_PolyDG, Heterodimer}, which employs the polygonal discontinuous Galerkin (PolyDG) method. This method offers distinct advantages due to its innate suitability for higher-order approximation and its versatility in mesh creation. Notably, it demonstrates exceptional effectiveness in approximating complex domains while upholding a superior level of accuracy, proving to be particularly valuable when addressing the intricate geometry of the brain \cite{polydg1}. Indeed, the PolyDG method empowers us to adjust approximation parameters locally, such as the polynomial degree $p$ and the element diameter $h$. This capability enables us to handle meshes featuring non-conforming elements with ease \cite{polydg3}. Furthermore, agglomeration-based strategies can be implemented, enhancing its versatility \cite{polydg4, polydg5}. This approach involves generating a coarser mesh to diminish the number of degrees of freedom in areas where it's unnecessary, thereby reducing the computational effort required. Such a process seamlessly integrates into the framework of polygonal and polyhedral grids due to their flexible element shape definitions \cite{polydg4, polydg5}.
\par
To demonstrate the practical capabilities of the proposed model we consider its application in the modeling of the onset of Alzheimer's disease.
This particular illness stands as one of the predominant forms of dementia, accounting for approximately $60\%$-$70\%$ of its cases \cite{WHO}. Alzheimer's disease is characterized by the agglomeration of protein fragments, specifically $\beta$-amyloid, forming extracellular neuritic plaques and twisted strands of $\tau$-protein, leading to intraneuronal neurofibrillary tangles \cite{ActaNeuro}. Observations indicate that the accumulation of this pathological material begins up to twenty years before the onset of the first symptoms of the disease \cite{ActaNeuro, Abetarate}. Our model can then be modified, albeit with some simplification, to depict this specific phenomenon or other neurodegenerative diseases with a similar origin.
\\
\par
The remaining part of the paper is organized as follows. In Section \ref{sec:MathModel2}, we introduce and present the mathematical model. In Section \ref{sec:PolyDG}, we outline the most important features of the PolyDG method for space discretization and introduce useful definitions necessary to our analysis. In Section \ref{sec:SD}, we derive the semi-discrete formulation of our problem by applying the PolyDG discretization on space. In Section \ref{sec:DisModel} we introduce the discretization in time and the two possible treatments of the nonlinear terms. In Section \ref{sec:Res}, we report and analyze the results of two convergence tests and simulations on three-dimensional domains. In Section \ref{sec: Simulations on brain}, we present simulations performed applying our model to the atrophy of the brain induced by Alzheimer's disease. To do so, we consider both the case in which we assume the atrophy to be an infinite process and the case of finite deformation applying a nonlinear elasticity equation. In Section \ref{sec:conclusion}, we present our conclusions and discuss future developments.
\section{Multiphysics coupled model of tissue atrophy and pathogen diffusion}\label{sec:MathModel2}
In this section, we construct a mathematical model describing the dynamics of pathogen diffusion and reaction coupled with a morpho-elastic description of tissue atrophy. Specifically, we adopt the FK equation to model the spreading of the species concentration and couple it with a description of the loss of the tissue mass using continuum mechanics. The strong formulation of the coupled model reads
\begin{equation}\label{Model}
    \begin{dcases}
        J\frac{\partial c}{\partial t} = \nabla_X \cdot (J\mathbf{F}^{-1}\mathbf{D}\mathbf{F}^{-T}\nabla_X c) + \alpha Jc(1-c) + Jf_c(\mathbf{x},t) & \text{in }\Omega \times (0,T],\\
        \dot{g} = \frac{1}{\tau}(g+1)\left(1-\frac{1}{\beta}(g+1)\right) & \text{in }\Omega\times (0,T],\\
        -\nabla_{X}\cdot\mathbf{P} + \mathbf{f_u} = \mathbf{0} & \text{in }\Omega \times (0,T],\\
        (\mathbf{D}\nabla c)\cdot \mathbf{n} = 0 & \text{on } \Gamma_{N}^c \times (0, T], \\
        c = c_D & \text{on } \Gamma_{D}^c \times (0, T], \\
        \mathbf{P}(\mathbf{u})\mathbf{n} = \mathbf{h_u} & \text{on } \Gamma_{N}^\mathbf{u} \times (0, T],\\
        \mathbf{u} = \mathbf{u_D} & \text{on } \Gamma_{D}^\mathbf{u} \times (0, T],\\
        c(\mathbf{x}, 0) = c_0 & \text{in }\Omega,\\
        g(\mathbf{x}, 0) = g_0 & \text{in }\Omega,\\ 
    \end{dcases}
\end{equation}
and is derived as follows.
\\
\subsection{Modeling pathogen dynamics}\label{subsec:ConModel}
We characterize the spreading of the pathogen by the FK equation \cite{Fisher,Kolmogorov}, which is frequently used in literature to model the propagation of a favored gene in population dynamics, \cite{physicsbasedmodelAlz}. We introduce the relative concentration of the pathogen as $c = c(\mathbf{x}, t):\Omega_t \times [0, T] \to \mathbb{R}$. 
The FK equation for the relative concentration $c=c(\mathbf{x}, t)$ can be formulated as follows:
\begin{equation}
\label{Fisher}
    \begin{cases}
        \frac{\partial c}{\partial t} = \nabla \cdot (\mathbf{D}\nabla c) + \alpha c(1-c) + f_c(\mathbf{x}, t) & \text{in }\Omega_t \times (0,T],\\
        (\mathbf{D}\nabla c)\cdot \mathbf{n} = 0 & \text{on } \Gamma_{Nt}^c \times (0, T], \\
        c = c_D & \text{on } \Gamma_{Dt}^c \times (0, T], \\
        c(\mathbf{x}, 0) = c_0(\mathbf{x}) & \text{in } \Omega_t.\\
    \end{cases}
\end{equation}
Here $f_c(\mathbf{x},t)$ is the forcing term, which models the external addition/removal of mass \cite{FK_equation_and_PolyDG}, $\mathbf{n}$ is the normal unit vector to the Neumann boundary, $c_D$ is the Dirichlet boundary condition, and $c_0(\mathbf{x})$ is the initial condition, which gives us the value and distribution of the concentration of the pathogen at the initial time $t=0$. The diffusion tensor $\mathbf{D}(\mathbf{x})$ describes the directions and velocity of the pathogen spreading in the tissue. We assume the diffusion tensor to be symmetric and positive definite. In equation \eqref{Fisher}, $\alpha = \alpha(\mathbf{x})$ is the reaction coefficient, modeling misfolding, clearance, and aggregation of the pathogens.
\par
Appropriate boundary conditions complement the FK equation. Specifically, we introduce a partition of the boundary $\partial\Omega_t$ into two subsets, denoted  by $\Gamma_{Nt}^c$ and $\Gamma_{Dt}^c$, so that $\partial \Omega_t = \Gamma_{Nt}^c \cup \Gamma_{Dt}^c$. We impose the Dirichlet boundary condition on $\Gamma_{Dt}^c$ so that $c =  c_D$. On $\Gamma_{Nt}^c$, we impose a homogeneous Neumann boundary condition that indicates the absence of flux of pathogens across the boundary.
\par
\begin{rmk} Assuming $f_c(\mathbf{x},t)=0$, for all $\mathbf{x}\in \Omega_t, t \in [0,T]$, $0\leq c_0(\mathbf{x})\leq 1$ and $c_0(\mathbf{x})\in H^1(\Omega_t)$ and taking homogeneous Neumann boundary conditions on the all boundary $\partial\Omega_t$, it is possible to prove that there exists a unique solution of \eqref{Fisher} such that $0\leq c(\mathbf{x}, t) \leq 1$ $\forall \mathbf{x}\in \Omega$, $\forall t > 0$ a.e.. In particular this means that, starting from a positive concentration, the solution propagates towards a stable equilibrium $c(\mathbf{x}, t) = 1$ with $t \to +\infty$ \cite{PDE}.
\end{rmk}
\par
Since we will consider an atrophy process we will consider a reference configuration $\Omega$, a current configuration $\Omega_t$ and a deformation $\phi$. I particular, equation~\eqref{Fisher} is written in the current configuration $\Omega_t$. However, when considering finite deformations, we need to solve~\eqref{Fisher} in the reference configuration $\Omega$. Then, the FK equation can be rewritten in a Lagrangian setting as
\begin{equation}\label{eq:FisherRef}
    J\frac{\partial c}{\partial t} = \nabla_X \cdot (J\mathbf{F}^{-1}\mathbf{D}\mathbf{F}^{-T}\nabla_X c) + \alpha Jc(1-c),
\end{equation}
where $\mathbf{F}$ is the deformation gradient, $J$ is the determinant of the deformation gradient and $\nabla_X$ the gradient operator in the reference system of coordinates $\mathbf{X}$.
The complete derivation of  \eqref{eq:FisherRef} can be found in Appendix \ref{AppB}, while the definition of $\Omega, \Omega_t$ and $\mathbf{F}$ will be introduced in the following section \eqref{GrModel}.
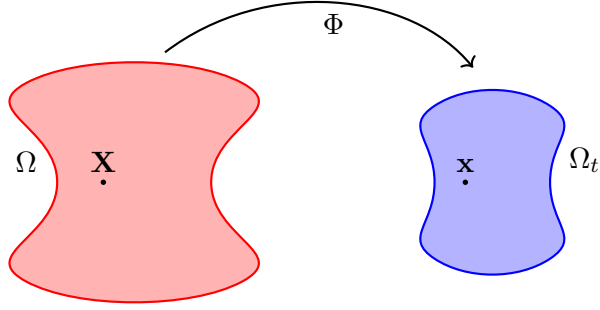
\begin{figure}[t!]
        \centering
         \resizebox{.5\textwidth}{!}{\begin{tikzpicture}
\fill [red, opacity = 0.3] plot [smooth cycle, tension = 1] coordinates {(2.4,0) (2,1.3) (4.8,1.3) (4.4,0) (4.8,-1.3) (2,-1.3)};
\draw [red, thick] plot [smooth cycle, tension = 1] coordinates {(2.4,0) (2,1.3) (4.8,1.3) (4.4,0) (4.8,-1.3) (2,-1.3)};
\fill (3, 0) circle (1pt) node[anchor=south] {$\mathbf{X}$};
\fill (2, 0) node[anchor=south] {$\Omega$};

\fill [blue, opacity = 0.3] plot [smooth cycle, tension = 1] coordinates {(7.3,0) (7.3,1) (8.8,1) (8.8,0) (8.8,-1) (7.3,-1)};
\draw [blue, thick] plot [smooth cycle, tension = 1] coordinates {(7.3,0) (7.3,1) (8.8,1) (8.8,0) (8.8,-1) (7.3,-1)};
\fill (7.7, 0) circle (1pt) node[anchor=south] {$\mathbf{x}$};
\fill (9.25, 0) node[anchor=south] {$\Omega_t$};

  \draw[->, thick]        (3.8,1.7) ..controls (5.1, 2.7) and (7, 2.5)  .. (7.8,1.5) node[midway, below] {$\Phi$};

\end{tikzpicture}}
    \caption{Graphic scheme of the reference domain $\Omega$ and the current domain $\Omega_t$ with the deformation $\Phi$.}
    \label{fig:AtrophyScheme}
\end{figure}

\subsection{Morpho-elasticity of tissue atrophy}\label{GrModel}

In this section, we introduce the morpho-elastic model for tissue atrophy, assuming that the size of the tissue undergoing a pathological loss of mass is much larger than the characteristic size of a cell. Therefore, we describe the tissue as a continuum elastic body, where an active mass modulation induces the deformation. 
\par
We assume that the tissue occupies a domain $\Omega\subset\mathbb{R}^d$ at $t=0$, which is assumed as the reference configuration, with $d=2,\,3$. The function that maps the reference domain to the current configuration $\Omega_t$ at time $t$ is the deformation $\bm{\phi}: \Omega\times (0,\,T] \rightarrow \mathbb{R}^d$.
We also introduce the displacement field $\mathbf{u}$, defined as the function $\mathbf{u}(t,\,\mathbf{X}) = \bm{\phi}(t,\,\mathbf{X})-\mathbf{X}$.
By $\mathbf{X}\in\Omega$, we indicate the generic point in the reference configuration. In contrast, we indicate with $\mathbf{x}$ the corresponding point in the current configuration at time $t$ so that $\mathbf{x} = \bm{\phi}(\mathbf{X},\,t)$. The deformation of the body can be described through the deformation gradient $\mathbf{F} = \nabla_\mathbf{X} \bm{\phi}$, where $\nabla_\mathbf{X}$ denotes the gradient operator with respect to the referential coordinates.

To model the mass reduction, we exploit the framework of morpho-elasticity \cite{Goriely_Math_Mec}. We consider a multiplicative decomposition of the deformation gradient \cite{rodriguez1994stress} so that
\begin{equation}
\label{multidec}
\mathbf{F} = \mathbf{F}_\text{e} \mathbf{G},
\end{equation}
where $\mathbf{F}_\text{e}$ represents the local elastic distortion of the material, while $\mathbf{G}$ describes the local inelastic distortion due to the growth or atrophy of the elastic body.
We next describe the mechanics of tissue atrophy and discuss the constitutive assumptions. Usually, mass loss takes place on a much longer time scale than the elastic deformations' time scale. Therefore, we can assume quasi-static deformations and neglect inertial effects, therefore balance of the linear momentum reads
\begin{equation}
\label{Elasticity}
    \begin{cases}
        - \nabla_\mathbf{X} \cdot \mathbf{P} + \mathbf{f_u} =\boldsymbol{0}& \text{in }\Omega \times (0,T],\\
        \mathbf{P}(\mathbf{u})\mathbf{n} = \mathbf{h_u} & \text{on } \Gamma_{N}^\mathbf{u} \times (0, T],\\
        \mathbf{u} = \mathbf{u_D} & \text{on } \Gamma_{D}^\mathbf{u} \times (0, T].
    \end{cases}
\end{equation}
where $\mathbf{P}$ is the first Piola-Kirchhoff stress tensor, $\nabla_\mathbf{X} \cdot$ is the divergence operator in the reference configuration and $\mathbf{f_u}$ is density of body forces. Additionally, we define $\Gamma^{\mathbf{u}}_{N}$ as the boundary region where we impose a traction $\mathbf{h_u}$, while $\Gamma^{\mathbf{u}}_{D}$ represents the portion of the boundary where we impose Dirichlet boundary conditions $\mathbf{u_D}$. As usual, we assume  $\Gamma^{\mathbf{u}}_{D} \cup \Gamma^{\mathbf{u}}_{N} = \partial\Omega$ and $\Gamma^{\mathbf{u}}_{D} \cap \Gamma^{\mathbf{u}}_{N} = \emptyset$.

In what follows, we describe the tissue as a hyperelastic material, i.e. we postulate the existence of a strain energy density $\Psi(\mathbf{X},\,\mathbf{F})$. We will omit the explicit dependence on $\mathbf{X}$ and $\mathbf{F}$ whenever convenient. Standard thermodynamic arguments \cite{gurtin2010mechanics} allows us to write
\[
\mathbf{P} = \frac{\partial\Psi}{\partial\mathbf{F}}.
\]
By following the standard theory of morpho-elasticity \cite{rodriguez1994stress,Goriely_Math_Mec}, the strain energy density $\Psi (\mathbf{X},\,\mathbf{F})$ of the material can be written as
\begin{equation}\label{Psi}
    \Psi(\mathbf{X},\,\mathbf{F},\,t) = \left(\det\mathbf{G}(\mathbf{X},\,t)\right)\Psi_0 (\mathbf{F}\mathbf{G}^{-1}(\mathbf{X},\,t)),
\end{equation}
where $\Psi_0$ represents the strain energy density of the material in its relaxed state. In this work we will assume to have an isotropic mass loss, i.e.
\begin{equation} \label{eq: G def}
    \mathbf{G} = (1 + g)\mathbf{I}.
\end{equation}
The expression of the Piola-Kirchhoff stress tensor becomes
\begin{equation}\label{Piola}
    \mathbf{P} = \det(\mathbf{G})\frac{\partial \Psi_0}{\partial \mathbf{F_e}}\mathbf{G}^{-T}.
\end{equation}
In what follows, we specialize the theoretical framework to the case of small elastic deformations. 

\subsubsection{Linearization of the mechanics of tissue atrophy}\label{sec:linear elasticity}
Let $\varepsilon$ be $\max_{\mathbf{X} \in \Omega}\max_{t\in (0,\,T]}\| \mathbf{u}(\mathbf{X},\,t) \|$, being $\| \cdot \|$ the Euclidean norm. The vector field $\mathbf{u_1}$ is the normalized counterpart of $\mathbf{u}$, i.e. $\mathbf{u}(\mathbf{X},\,t) = \varepsilon \mathbf{u_1}(\mathbf{X},\,t)$, where $\varepsilon$ is a real positive number assumed small. Therefore, we can take an series expansion of $\mathbf{G}$ \eqref{eq: G def} in $\varepsilon$, assuming $\mathbf{G}$ to be a small perturbation of the identity, i.e.
\begin{equation}\label{G definition}
    \mathbf{G} = \mathbf{I} + \varepsilon \mathbf{G_1} + o(\varepsilon).
\end{equation}
This definition of $\mathbf{G}$ is aligned with \eqref{eq: G def}, taking $\varepsilon \mathbf{G_1} = g\mathbf{I}$, and choosing the parameter $\gamma$ \eqref{g_bar equation}, appropriately small.
Using the expansion\eqref{G definition} of $\mathbf{G}$ we obtain:
\begin{equation}
    \det \mathbf{G} = \det(\mathbf{I} + \varepsilon \mathbf{G_1} + o(\varepsilon)) = 1 + \varepsilon \tr \mathbf{G_1} + o(\varepsilon),
\end{equation}
\begin{equation}
    \mathbf{G}^{-1} = \mathbf{I} - \varepsilon \mathbf{G_1} + o(\varepsilon),
\end{equation}
\begin{equation}
    \mathbf{F}\mathbf{G}^{-1} = (\mathbf{I} + \varepsilon \nabla \mathbf{u_1})(\mathbf{I} - \varepsilon\mathbf{G_1} + o(\varepsilon) = \mathbf{I} + \varepsilon(\nabla\mathbf{u_1} - \mathbf{G_1}) + o(\varepsilon).
\end{equation}
Performing a Taylor expansion on \eqref{Piola}, the Piola-Kirchoff stress tensor can be approximated as:
\begin{equation}\label{P order 1}
    \mathbf{P} = \varepsilon\mathbb{C}:(\nabla\mathbf{u_1} - \mathbf{G_1}) + o(\varepsilon)= \mathbb{C}:(\nabla\mathbf{u} - g\mathbf{I}) = \mathbb{C}:\mathbf{E_e}
\end{equation}
where we neglect the higher order terms in $\varepsilon$ and we take $\mathbf{E_e} = \frac{\nabla\mathbf{u} - g\mathbf{I} + (\nabla\mathbf{u} - g\mathbf{I})^T}{2}$.
For an isotropic material, it is well-known that $\mathbf{P} = \mathbb{C}:\mathbf{E_e} = 2\mu\mathbf{E_e} + \lambda(\tr(\mathbf{E_e}))\mathbf{I}$, which gives
\begin{equation}\label{Piola_Kirchoff}
         \mathbf{P} = 2\mu\left(\frac{\nabla \mathbf{u} + \nabla\mathbf{u}^T}{2}\right) + \lambda \nabla \cdot \mathbf{u}\mathbf{I} - 2\mu\left(\frac{g\mathbf{I} + (g\mathbf{I})^T}{2}\right) - \lambda \tr\left(\frac{g\mathbf{I} + (g\mathbf{I})^T}{2}\right)\mathbf{I}.
\end{equation}
Under this assumption the Piola-Kirchhoff tensor simplifies into
\begin{equation}\label{eq: Piola def}
    \mathbf{P} = 2\mu\left(\frac{\nabla \mathbf{u} + \nabla\mathbf{u}^T}{2}\right) + \lambda \nabla \cdot \mathbf{u}\mathbf{I} -(2\mu + d\lambda)g \mathbf{I},
\end{equation}
where we recall $d$ is the space dimension.

Defining $\mathbf{E}(\mathbf{u}) = \frac{1}{2}(\nabla \mathbf{u} + \nabla\mathbf{u}^T)$, the symmetric part of the gradient of the displacement, we can write the elasticity equation as:
\begin{equation} \label{Pform}
    \begin{split}
        -\nabla\cdot\mathbf{P} &= -2\nabla\cdot(\mu\mathbf{E}(\mathbf{u})) -\nabla\cdot(\lambda \nabla \cdot \mathbf{u}\mathbf{I}) + \nabla\cdot((2\mu + d\lambda)g \mathbf{I}) \\ &=-2\nabla\cdot(\mu\mathbf{E}(\mathbf{u})) -\nabla(\lambda \nabla \cdot \mathbf{u}) + \nabla((2\mu + d\lambda)g) = \mathbf{f_u}.
    \end{split}
\end{equation}
We remark that the local change of volume due to atrophy is given by $\det \mathbf{G}$. In the linear elastic setting introduced in this section we get $\det \mathbf{G} = 1 + \tr g\mathbf{I} + o (\varepsilon) = 1+3g +o(\varepsilon)$. Thus, $3g$ serves as an indicator of the local volume change.

In a linear deformation regime, we can neglect the difference between the current domain $\Omega_t$ and the reference domain $\Omega$, therefore, we assume $\Omega \simeq \Omega_t$, $\partial \Omega \simeq \partial\Omega_t$.

\subsubsection{Evolution law for the inelastic tensor $\mathbf{G}$}\label{sec:Atrophy}
The evolution law for the inelastic tensor $\mathbf{G}$ introduced in equation \eqref{multidec} must be constitutively provided. We assume that the atrophy of the tissue is isotropic, with $\mathbf{G} = (1+g) \mathbf{I}$. The scalar variable $g = g(\mathbf{X},\,t)$ accounts for the local mass reduction triggered by the local concentration of the pathogen. Specifically, the mass density at time $t$ in the reference configuration is given by
\begin{equation}\label{eq:log. ode}
    \rho_t(\mathbf{X},\,t) = \rho_0(\mathbf{X}) \det\mathbf{G} (\mathbf{X},\,t),
\end{equation}
where $\rho_0$ is the initial mass density in the reference configuration. In this model, we assume that no mass change occurs if the pathogen's concentration $c(\mathbf{x}, t)$ is below a critical threshold $c_\text{cr}$. Tissue atrophy takes place when the concentration $c$ is above such a critical threshold. To mimic this behavior, we introduce the following logistic-type equation for the evolution of $g(t)$
\begin{equation}\label{logistic equation}
    \left\{
    \begin{aligned}
        &\dot{g} = \frac{1}{\widehat{\tau}}(1+g)\left(1-\dfrac{1+g}{\beta}\right) && \text{in }(0,T],\\
        &g(\mathbf{x}, 0) = 0 && \text{in }\Omega,\\ 
        \end{aligned}
    \right.
\end{equation}
where 
\begin{equation} \label{g_bar equation}
    \beta = \begin{cases}
        1 & \text{if $c\le c_{cr}$},\\
        1 - \gamma\frac{c - c_{cr}}{1 - c_{cr}} & \text{if $c>c_{cr}$}.
    \end{cases}
\end{equation}
The parameter $\widehat{\tau}$ in equation \eqref{logistic equation} represents the characteristic time of tissue atrophy, while $1-\gamma$ in equation \eqref{g_bar equation}, with $0<\gamma<1$, represents the minimal value of $\beta$, obtained for $c(\mathbf{x}, t)=1$.
\begin{rmk}
    (Stability analysis) We compute the equilibria of equation \eqref{logistic equation} by imposing the right-hand side of the equation equal to $0$, and find $g+1= \{ \beta, 0 \}$. Now we observe that by evaluating the derivative $\frac{\partial G(g)}{\partial g}$ in $g = \beta - 1$ we obtain:
\begin{equation}
    \frac{1}{\widehat{\tau}}\left(1-\frac{2}{\beta}(1+g)\right)\big|_{g = \beta-1} = -\frac{1}{\widehat{\tau}} < 0,
\end{equation}
which implies that $g = \beta - 1$ is a stable equilibrium for the logistic equation.
Furthermore, we observe that, when $c \leq c_{cr}$, the stable equilibrium $g = \beta - 1 = 0$. In fact, we do not have atrophy caused by the pathogen. When, instead, $c>c_{cr}$, the stable equilibria become $g = \beta - 1 = - \gamma \frac{c - c_{cr}}{1 - c_{cr}}$, which varies with the increase of the concentration of the pathogen. In particular, when $c = 1$, $g = \beta - 1 = -\gamma$, meaning that $\gamma$ represents the absolute value of the maximum shrinkage that we can obtain.
\par
Now, coupling together equations \eqref{eq:FisherRef}, \eqref{logistic equation} and \eqref{Elasticity} we recover the strong formulation of the coupled problem~\eqref{Model}.
\end{rmk}

\subsection{Weak formulation}
Now we derive the variational formulation of problem \eqref{Model}. To begin with, we introduce functional spaces tailored for the solutions of the equations. To address the solution of the parabolic differential FK equation, we define the spaces: $W = H^1(\Omega)$, $W_0 = H^1_{\Gamma_D^c}(\Omega) = \{w \in H^1(\Omega) : w|_{\Gamma_D^c} = 0\}$ and $W_D = \{w \in H^1(\Omega) : w|_{\Gamma_D^c} = c_D\}$. Turning to the solution $g$ of the ODE~\eqref{logistic equation}, we introduce the space $Q = H^1(\Omega)$. Additionally, for the solution of the elasticity equation~\eqref{Elasticity}, we define $\mathbf{V} = \mathbf{H}^1(\Omega, \mathbb{R}^d)$, $\mathbf{V}_0 = \mathbf{H}^1_{\Gamma_D^\mathbf{u}}(\Omega, \mathbb{R}^d):=\{\mathbf{v} \in \mathbf{H}^1(\Omega, \mathbb{R}^d) : \mathbf{v}|_{\Gamma_D^\mathbf{u}} = 0\}$, and $\mathbf{V}_D :=\{\mathbf{v} \in \mathbf{H}^1(\Omega, \mathbb{R}^d) : \mathbf{v}|_{\Gamma_D^\mathbf{u}} = \mathbf{u_D}\}$.
Furthermore, we utilize the conventional definition of the $L^2(\Omega)$ scalar product, represented as $(\cdot, \cdot)_{\Omega}$ with the associated norm denoted by $||\cdot||_{\Omega}$. This definition extends componentwise for vector-valued and tensor-valued functions \cite{PartialDiffEq}.
The weak formulation of the problem described in System \eqref{Model} becomes:\\
$\forall t \in (0,T]$ find $c(\mathbf{x}, t) \in W_D$, $g(\mathbf{x},t) \in Q$ and $\mathbf{u}(\mathbf{x})\in\mathbf{V}_0$ such that:
\begin{equation}\label{WeakModel}
    \begin{cases}
        \left(J\frac{\partial c(\mathbf{x},t)}{\partial t}, w \right)_\Omega + a_c(c(\mathbf{x},t), w) - r_L^c(c(\mathbf{x},t), w) + r_N^c(c(\mathbf{x},t), c(\mathbf{x},t), w) = F_c(w) & \forall w \in W_0,\\
        \left(\dot{g}(\mathbf{x}, t),p\right)_\Omega = r_L^g(g,p) - r_N^g(g,g,p) + F_g(p) & \forall p \in Q,\\
        a_E(\mathbf{u},\mathbf{v}) = F_E(\mathbf{v}) & \forall \mathbf{v} \in \mathbf{V}_0,\\
        c(\mathbf{x},0) = c_0(\mathbf{x}) & \text{in } \Omega,\\
        g(\mathbf{x},0) = 0 & \text{in } \Omega,
    \end{cases}
\end{equation}
where:
\begin{align*}
    &a_c(c,w) = (J\mathbf{F}^{-1}\mathbf{D}\mathbf{F}^{-T}\nabla c, \nabla w )_{\Omega}, &&r_L^c(c,w) = \left(J\alpha c,w\right)_{\Omega},\\
    &r_N^c(v, c, w) = \left(J\alpha (vc), w\right)_{\Omega},
    && F_c(w) = \left(Jf_c, w\right)_{\Omega},\\
    &r_L^g(g,p) = \left(\frac{1}{\widehat{\tau}}\left(1-\frac{2}{\beta(c)}\right)g,p \right)_\Omega, &&r_N^g(g,p,q) = \left(\frac{1}{\widehat{\tau} \beta(c)}g,p,q\right)_\Omega,\\
    &F_g(p) = \frac{1}{\widehat{\tau}} \left( \left(1 - \frac{1}{\beta(c)} \right), p \right)_\Omega, &&a_E(\mathbf{u},\mathbf{v}) = \left(\mathbf{P}(\mathbf{u}),\nabla_{X}\mathbf{v}\right)_{\Omega},\\ &F_E(\mathbf{v}) = (\mathbf{f_u}, \mathbf{v})_\Omega + (\mathbf{h_u},\mathbf{v})_{\Gamma_N^\mathbf{u}}. &&\\ 
\end{align*}
For all $c, w, v \in W$, $g, p, q \in Q$ and $\mathbf{u}, \mathbf{v} \in \mathbf{V}$.

\section{Polygonal Discontinuous Galerkin formulation}\label{sec:PolyDG}
Let $\mathscr{T}_h$ a partition of the domain $\Omega$ using polygonal or polyhedral elements $K \in \mathscr{T}_h$. Here, $|K|$ denotes the measure of each element, $h_K$ represents the diameter of each element, and $h = \max_{K \in \mathscr{T}_h} h_K < 1$. The interfaces of each element are the ($d-1$)-dimensional intersections of adjacent facets. We distinguish two cases:
\begin{itemize}
    \item case $d = 3$, in which the interfaces consists in triangles, the set of which we denote by $\mathscr{F}^h$;
    \item case $d = 2$, in which the interfaces are line segments, the set of which we denote by $\mathscr{F}^h$. 
\end{itemize} Specifically, the set of interfaces $\mathscr{F}_h$ comprises the union of boundary faces $\mathscr{F}^B_h$, lying on the boundary, and all interior faces $\mathscr{F}^I_h$.  We can further categorize the set of boundary faces $\mathscr{F}^B_h$ into interfaces where Dirichlet conditions are applied, $\mathscr{F}^D_h$, and interfaces where Neumann boundary conditions are applied, $\mathscr{F}^N_h$. We also assume that $\mathscr{T}_h$ is aligned with $\Gamma_D$ and $\Gamma_N$, implying that any element in $\mathscr{F}^B_h$ is contained within either $\Gamma_D$ or $\Gamma_N$. For additional assumptions on polytopal meshes.
\par
We introduce the trace operators on the interior faces $\mathscr{F}^I_h$. We also employ the notation ($\pm$) to signify traces of functions on $F \in \mathscr{F}_h^I$ within the interior of $K^{\pm}$ for a generic function. For a scalar-valued function $q$, a vector-valued function $\mathbf{v}$ and a tensor-valued function $\bm{\tau}$ we define:
\begin{itemize}
    \item $\averagel q\averager = \frac{1}{2}(q^+ + q^-)$, $\averagel \mathbf{v}\averager = \frac{1}{2}(\mathbf{v}^+ + \mathbf{v}^-)$ and $\averagel \bm{\tau}\averager = \frac{1}{2}(\bm{\tau}^+ + \bm{\tau}^-)$
    \item $\jumpl \cdot \jumpr$ on $F \in \mathscr{F}^I_h$: $\jumpl q \jumpr = q^+\mathbf{n}^+ + q^-\mathbf{n}^-$, $\jumpl\mathbf{v}\jumpr = \mathbf{v}^+\cdot\mathbf{n}^+ + \mathbf{v}^-\cdot\mathbf{n}^-$, and $\jumpl \bm{\tau} \jumpr = \bm{\tau}^+\mathbf{n}^+ + \bm{\tau}^-\mathbf{n}^-$
    \item $\jjumpl \mathbf{v} \jjumpr = \frac{1}{2}(\mathbf{v}^+\otimes\mathbf{n}^+ + \mathbf{n}^+\otimes\mathbf{v}^+) + \frac{1}{2}(\mathbf{v}^-\otimes\mathbf{n}^- + \mathbf{n}^-\otimes\mathbf{v}^-)$
\end{itemize}
On $F \in \mathscr{F}^D_h$, we set the following trace operators for the test functions as $\averagel q\averager = q$, $\averagel\mathbf{v}\averager = \mathbf{v}, \averagel\bm{\tau}\averager = \bm{\tau}$ and $\jjumpl\mathbf{v}\jjumpr = \frac{1}{2}(\mathbf{v}\otimes\mathbf{n} + \mathbf{n}\otimes\mathbf{v})$. Additionally, for the trial functions we define the traces operator on the faces of the Dirichlet boundary as $\jumpl p \jumpr = (p-h_D)\mathbf{n}$, $\jumpl\mathbf{u}\jumpr = (\mathbf{u}-\mathbf{h_D})\cdot\mathbf{n}, \ \jumpl\bm{\tau}\jumpr = (\bm{\tau} - \bm{\gamma}_D) \mathbf{n}$ and $\jjumpl\mathbf{v}\jjumpr = \frac{1}{2}((\mathbf{v} - \mathbf{h_D})\otimes\mathbf{n} + \mathbf{n}\otimes(\mathbf{v} - \mathbf{h_D}))$, with $\otimes$ defined as $\mathbf{a}\otimes\mathbf{b} \coloneqq \mathbf{a}\mathbf{b}^T$ meaning $(\mathbf{a} \otimes \mathbf{b})_{ij} \coloneqq a_ib_j \mathbf{e}_i \otimes \mathbf{e}_j$.
\par
We define the penalty functions $\eta: \mathscr{F}_h \to \mathbb{R}$ and $\xi:\mathscr{F}_h \to \mathbb{R}_+$ defined faced-wise as:
\begin{flalign*}
    \eta = \eta_0 \left\{\begin{array}{ll}
        \max\{\{d^K\}_H, \{\alpha\}\}\frac{p^2}{\{h\}_H}, & on \ F \in\mathscr{F}^I_h, \\
        \max\{d^K,\alpha \}\frac{p^2}{h}, & on \ F \in \mathscr{F}^D_h,
    \end{array}\right. && \xi = \xi_0 \left\{\begin{array}{ll}
        \{\tilde{\mathbb{C}}_E^K\}_H\frac{p^2}{\{h\}_H}, & on \ F \in\mathscr{F}^I_h, \\
        \tilde{\mathbb{C}}_E^K\frac{p^2}{h}, & on \ F \in \mathscr{F}^D_h,
    \end{array}\right.
\end{flalign*}
where $\{\cdot\}_H$ denotes the harmonic average operator on $F \in \mathscr{F}^I_h$, $\eta_0$ and $\xi_0$ are constants chosen sufficiently large to guarantee the stability of the methods, $d^K$ and $\tilde{\mathbb{C}}_E^K$ are defined as $d^K := \left\|\sqrt{\Tilde{\mathbf{D}}|_K}\right\|^2$ and  $\tilde{\mathbb{C}}_E^K \coloneqq \left\|\sqrt{\mathbb{C}_E|_K}\right\|^2$ for any $K \in \mathscr{T}_h$.\\
We define: $W_h^{DG} \coloneqq \{ w \in L^2(\Omega) : w_{|_K} \in \mathbb{P}_r(K) \ \forall K \in \mathscr{T}_h \} $, $Q_h^{DG} \coloneqq \{ q \in L^2(\Omega) : q_{|_K} \in \mathbb{P}_q(K) \ \forall K \in \mathscr{T}_h \}$ and $\mathbf{V}_h^{DG} \coloneqq \{ \mathbf{v} \in L^2(\Omega; \mathbb{R}^d) : \mathbf{v}_{|_K} \in [\mathbb{P}_p(K)]^d \ \forall K \in \mathscr{T}_h \}$. 

\section{PolyDG semi-discrete formulation} \label{sec:SD}
In this section we derive the PolyDG semi-discrete formulation of the weak problem \eqref{WeakModel}.
Setting $\int_\mathscr{F_h} = \sum_{F\in\mathscr{F_h}} \int_F$, we define the following bilinear forms:\\
\begin{itemize}
    \item $\mathscr{A}_c:W^{DG}_h \times W^{DG}_h \to \mathbb{R}$:
    \begin{equation}
        \begin{split}
            \mathscr{A}_c(c_h, w_h) \coloneqq & \int_\Omega \mathbf{D}\nabla_h c_h \cdot \nabla_h w_h + \int_{\mathscr{F}^I_h \cup \mathscr{F}^D_h} \eta \jumpl c_h \jumpr\cdot \jumpl w_h \jumpr \mathrm{d}\sigma \\ 
            & -\int_{\mathscr{F}^I_h \cup \mathscr{F}^D_h} (\averagel\mathbf{D}\nabla_h c\averager\cdot\jumpl w_h \jumpr+\jumpl c_h \jumpr\cdot\averagel\mathbf{D}\nabla_h w_h\averager)\mathrm{d}\sigma \quad \forall c_h,w_h \in W^{DG}_h
        \end{split}
    \end{equation}
    \item $\mathscr{A}_E: \mathbf{V}^{DG}_h \times \mathbf{V}^{DG}_h \to \mathbb{R}$:
    \begin{equation}
        \begin{split}
            \mathscr{A}_E(\mathbf{u_h}, \mathbf{v_h}) \coloneqq &\int_\Omega 2\mu\bm{\varepsilon}(\mathbf{u_h}):\bm{\varepsilon}(\mathbf{v_h}) + \lambda\nabla_h\cdot\mathbf{u_h}\nabla_h\cdot\mathbf{v_h} + \int_{\mathscr{F}^I_h \cup \mathscr{F}^D_h} \xi\jjumpl \mathbf{u_h}\jjumpr:\jjumpl\mathbf{v_h}\jjumpr \mathrm{d}\sigma\\
            & -\int_{\mathscr{F}^I_h \cup \mathscr{F}^D_h} ( \averagel \mathbf{P}_E(\mathbf{u_h}) \averager:\jjumpl\mathbf{v_h}\jjumpr + \jjumpl\mathbf{u_h}\jjumpr:\averagel\mathbf{P}_E(\mathbf{v_h})\averager)\mathrm{d}\sigma \quad \forall \mathbf{u_h},\mathbf{v_h}\in \mathbf{V}_h^{DG},
        \end{split}
    \end{equation}\\
    \item $\mathscr{B}_E: Q_h^{DG} \times \mathbf{V}^{DG}_h \to \mathbb{R}$:
    \begin{equation}
         \mathscr{B}_E(p_h, \mathbf{v_h}) = \int_\Omega (2\mu + d\lambda)p_h \nabla_h \cdot \mathbf{v_h} - \int_{\mathscr{F}^I_h \cup \mathscr{F}^D_h} (2\mu + d\lambda)\averagel p_h \mathbf{I}\averager\!:\!\jjumpl\mathbf{v_h}\jjumpr \mathrm{d}\sigma \quad \forall p_h \in Q_h^{DG}, \mathbf{v}_h \in \mathbf{V}_h^{DG}.
    \end{equation}
\end{itemize}
In the above definitions, $\nabla_h$ denotes the elementwise gradient \cite{NumModels}. Employing the aforementioned bilinear forms, we derive the semi-discrete formulation of model \eqref{WeakModel} as follows:
\par
\bigskip
For any $t \in (0,T]$ find $c_h(t) \in W^{DG}_h, g_h(t) \in Q_h^{DG}$ and $\mathbf{u}_h(t) \in \mathbf{V}_h^{DG}$ such that:
\begin{equation}\label{SDModel}
    \begin{cases}
        \left(\dot{c}_h, w_h \right)_\Omega + \mathscr{A}_c(c_h, w_h) - r_L(c_h, w_h) + r_N(c_h, c_h, w_h) = F_c(w_h) & \forall w_h \in W^{DG}_h,\\
        \left(\dot{g}_h, p_h\right)_\Omega = r_L^g(g_h, p_h)-r_N^g(g_h,g_h,p_h) + F_g(p_h) & \forall p_h \in Q_h^{DG},\\
        \mathscr{A}_E(\mathbf{u_h}, \mathbf{v_h}) - \mathscr{B}_E(g_h, \mathbf{v_h}) = F_E(\mathbf{v}_h) & \forall \mathbf{v}_h \in \mathbf{V}_h^{DG},\\
        c_{h}(0) = c_{0h}, \qquad g_h(0) = 0 & \mathrm{in} \, \ \Omega_h.\\
    \end{cases}
\end{equation}
In Equation \eqref{SDModel}, $c_{0h}$ is a suitable approximation of the initial conditions $c_0$ in the discrete space $W_h^{DG}$.

\subsection{Algebraic formulation}
\label{subsec:Algform}
Let $\{\phi_{j}\}^{N_c}_{j=1}$, $\{q_j\}_{j=1}^{N_g}$, and $\{\bm{\psi}_j\}_{j=1}^{N_\mathbf{u}}$ be suitable bases for the discrete spaces $W^{DG}_h$, $Q^{DG}_h$, and $\mathbf{V}^{DG}_h$, respectively. Then we can write:
\begin{equation*}
    c_h(t) = \sum_{j=1}^{N_c}C_n(t)\phi_j, \qquad g_h(t) = \sum_{j=1}^{N_g}g_n(t)q_j, \qquad \mathbf{u}_h = \sum_{j=1}^{N_\mathbf{u}}U_n(t)\bm{\psi}_j.
\end{equation*}
We denote $\mathbf{C} \coloneqq [C_n]_{n=1}^{N_c} \in \mathbb{R}^{N_c}$, $\mathbf{g} \coloneqq [G_n]_{n=1}^{N_g} \in \mathbb{R}^{N_g}$ and $\mathbf{U} \coloneqq [U_n]_{n=1}^{N_\mathbf{u}} \in \mathbb{R}^{dN_\mathbf{u}}$ the vector of the expansion coefficients. The algebraic form of 
~\eqref{SDModel} can be written as:
Find $\mathbf{C}(t)\in\mathbb{R}^{N_c}$, $\mathbf{g}(t)\in\mathbb{R}^{N_g}$ and $\mathbf{U}\in\mathbb{R}^{dN_\mathbf{u}}$ such that $\forall t \in (0,T]$ we have:
\begin{equation}\label{algebraicModel}
    \begin{cases}
        \mathbf{M}_c\dot{\mathbf{C}}(t)+\mathbf{A}_c\mathbf{C}(t)-\mathbf{M}_\alpha\mathbf{C}(t)+\widetilde{\mathbf{M}}_\alpha(\mathbf{C}(t))\mathbf{C}(t) = \mathbf{F}_c(t), & t \in (0,T],\\
         \mathbf{M}_g\dot{\mathbf{g}}(t)=\mathbf{M}_g\mathbf{g}(t)-\widetilde{\mathbf{M}}_g(\mathbf{g}(t))\mathbf{g}(t) + \mathbf{F}_g& t \in (0,T],\\
         \mathbf{K}_E\mathbf{U} - \mathbf{B}_g^T \mathbf{g} = \mathbf{F}_E,\\
        \mathbf{C}(0) = \mathbf{C}_0,\\
        \mathbf{g}(0) = \mathbf{g}_0.
    \end{cases}
\end{equation}
The reader can find the definition of the matrices in Appendix \ref{AppA}. 
\section{Fully discrete formulation}\label{sec:DisModel}
Now, let us present the fully discrete approximation of Equation~\eqref{Model}. We apply the Crank-Nicolson method to discretize temporal derivatives, and we consider a semi-implicit treatment of the nonlinear terms, ultimately outlining the complete discrete formulation of our problem.
To discretize the time evolution of \eqref{algebraicModel} we define a partition of the time interval $[0, T]$ into $N$ intervals: $0=t_0<t_1<...<t_N=T$ and we assume a time step $\Delta t = t_{n+1} - t_{n}, n=0,...,N-1$. 
Consequently, the fully discrete formulation of Problem \eqref{Model} becomes: Given $\mathbf{C}(0) = \mathbf{C}_0$ and $\mathbf{g}(0) = \mathbf{g}_0$, find $\mathbf{C}^{n+1}\simeq\mathbf{C}(t_{n+1})$, $\mathbf{g}^{n+1}\simeq\mathbf{g}(t_{n+1})
$ and $\mathbf{U}^{n+1}\simeq\mathbf{U}(t_{n+1}) \in \mathbb{R}^{dN_\mathbf{u}}$ for $n = 1,...,N-1$:
\begin{equation}\label{CNModel}
    \begin{cases}
        \mathbf{M}_c\mathbf{C}^{n+1}+\frac{\Delta t}{2}(\mathbf{A}_c - \mathbf{M}_\alpha)\mathbf{C}^{n+1}+\Delta t \widetilde{\mathbf{M}}_\alpha(\mathbf{C}^*)\mathbf{C}^{n+1,n} = \mathbf{M}\mathbf{C}^{n} - \frac{\Delta t}{2}(\mathbf{A} - \mathbf{M}_\alpha)\mathbf{C}^{n}+\frac{1}{2}(\mathbf{F}_c^{n+1}+\mathbf{F}_c^{n}),\\
        \mathbf{M}_g\mathbf{g}^{n+1}-\frac{\Delta t}{2}\mathbf{M}_\beta\mathbf{g}^{n+1}+\Delta t\widetilde{\mathbf{M}}_\beta(\mathbf{g}^*)\mathbf{g}^{n+1,n} = \mathbf{M}_g\mathbf{g}^{n}+\frac{\Delta t}{2}\mathbf{M}_\beta\mathbf{g}^{n} + \frac{\Delta t}{2}(\mathbf{F}_g^{n+1} + \mathbf{F}_g^{n}),\\
        \mathbf{K}_E\mathbf{U}^{n+1} - \mathbf{B}_g^T \mathbf{g}^{n+1} = \mathbf{F}_E^{n+1},
    \end{cases}
\end{equation}
where  $\mathbf{C}^{n+1,n} = \frac{1}{2}(\mathbf{C}^{n+1}+\mathbf{C}^{n})$ and $\mathbf{g}^{n+1,n} = \frac{1}{2}(\mathbf{g}^{n+1}+\mathbf{g}^n)$, and  $\mathbf{C}^* \coloneqq \frac{3}{2}\mathbf{C}^{n}-\frac{1}{2}\mathbf{C}^{n-1}$ and $\mathbf{g}^* = \frac{3}{2}\mathbf{g}^{n}-\frac{1}{2}\mathbf{g}^{n-1}$.
\par
\bigskip
In the end, we summarize the basic idea of the complete numerical solver in Algorithm \ref{alg:1}.
\begin{algorithm}
    \caption{Numerical Algorithm}\label{alg:1}
    \begin{algorithmic}
        \REQUIRE $c_{h0}$ and $g_{h0}$
        \WHILE{$t_n \le T$}
        \STATE Solve FK $\rightarrow c_h^n$ \ \ \ \eqref{Fisher}
        \IF{$c_h^n < c_{cr}$} 
            \STATE $\beta \gets 1$
        \ELSE
        \STATE $\beta \gets \beta(c_h^n)$
        \ENDIF
        \STATE Solve Logistic Law with parameter $\beta$ $\rightarrow g_h^n$ \ \ \ \eqref{eq:log. ode}
        \STATE Use $g_h^n$ to define the growth tensor and  $ \mathbf{P}(\mathbf{u}_h^n)$ \ \ \ \eqref{eq: Piola def}
        \STATE Solve Elasticity Equation $\rightarrow \mathbf{u}_h^n$ \ \ \ \eqref{Elasticity}
        \STATE $t_n \rightarrow t_{n+1}$
        \ENDWHILE
    \end{algorithmic}
\end{algorithm}
\section{Numerical results}\label{sec:Res}
In this section, we present the results of a numerical convergence test conducted to demonstrate the accuracy of the proposed method.
\par
The numerical simulations are done using the FEniCS finite element software \cite{AlnaesEtal2015}. This software allows us to perform numerical simulations on a 3-dimensional domain with a mesh of tetrahedral elements. We use the MUMPS solver to solve the FK and the logistic equation. In contrast, we employ the iterative solver GMRES with a SOR preconditioner for the elasticity equation.
\par
To perform a convergence test we denote the L$^2$ norm as $|| \cdot || := || \cdot ||_{L^2}$ on $\Omega$, and the $L^2$-norm on a set of faces $\mathscr{F}$: $||\cdot||_{\mathscr{F}}:= \bigl(\sum_{F \in \mathscr{F}} ||\cdot||_{L^2(F)}^2 \bigr)^{\frac{1}{2}}$ it is possible to define the DG-norms as: 
\begin{equation} 
    \|c\|_{DG} = \|\sqrt{\mathbf{D}}\nabla_h c\| + \|\sqrt{\eta}\jumpl c \jumpr\|_{\mathscr{F}^I_h \cup \mathscr{F}^D_h} \qquad \forall c \in H^1(\mathscr{F}_h)
\end{equation} 
\begin{equation}\label{uDGnorm}
    \|\mathbf{u}\|_{DG} = \|\sqrt{\mathbb{C}_E}[\bm{\varepsilon}_h(\mathbf{u})]\|+\|\sqrt{\xi}[\![\![\mathbf{u}]\!]\!]\|_{\mathscr{F}^I_h \cup \mathscr{F}^D_h} \qquad \forall \mathbf{u} \in \mathbf{H}^1(\mathscr{T}_h, \mathbb{R}^d)
\end{equation}
\subsection{Test case 1: convergence analysis of a three-dimensional test case}
\label{sec:test_case_1}
In this test case, we perform a convergence analysis employing the parameters shown in Table \ref{tab:test2_parameters}, notice that this set of parameters is consistent with the 
forthcoming application of Section~\eqref{sec: Simulations on brain} \cite{FK_equation_and_PolyDG, Poroelasticity_PolyDG}. 

\begin{table}[t]
\caption{Parameters used in the test cases 1 and 2.}\label{tab:test2_parameters}
\begin{tabular}{lrlc|lrlc}
\toprule
  \textbf{Parameters} & \multicolumn{2}{c}{\textbf{Values}} &
  \textbf{Reference} & \textbf{Parameters} & \multicolumn{2}{c}{\textbf{Values}} &
  \textbf{Reference} \\ 
\midrule
 \normalsize$d_{ext}$ & \normalsize$8.00$& $[\mathrm{mm^2/year}]$ & \normalsize\cite{biochemical_biomechanical_Kuhl} & \normalsize$\lambda$ &
 \normalsize $505$ & $[\mathrm{Pa}]$&
 \normalsize \cite{Poroelasticity_PolyDG}\\
 \normalsize$d_{axn}$ & \normalsize $0.00 $ &$[\mathrm{mm^2/year}]$ & \normalsize \cite{biochemical_biomechanical_Kuhl} & \normalsize $\mu$ & \normalsize $216$ &$[\mathrm{Pa}]$ &
 \normalsize \cite{Poroelasticity_PolyDG} \\
 \normalsize $\alpha$ &
 \normalsize $0.90$ & $[-]$ & \normalsize \cite{biochemical_biomechanical_Kuhl} & \normalsize $\widehat{\tau}$ &\normalsize $1$& $[\mathrm{year}]$ &
 \normalsize  \\
  \normalsize $\gamma$ & \normalsize $0.05$& $[-]$ & \normalsize \\ 

\bottomrule
\end{tabular}
\end{table}
We consider a cubic domain $\Omega = (0, 1)^3$ and a time interval $[0, 0.1]$ with a time step $\Delta t = 10^{-3}$. We impose as exact solution:
\begin{align*}
        &c(\mathbf{X},t) = \textrm{cos}(\pi X)\textrm{cos}(\pi Y)\textrm{cos}(\pi Z)e^{-t}, &&
        g(\mathbf{X},t) = 1.0/(2e^{t/\tau} - 1.0),\\
        &\mathbf{u}(\mathbf{X}) = \left[
    \begin{array}{c}
         -\textrm{cos}(2\pi X)\textrm{cos}(2\pi Y)  \\
         \textrm{sin}(2\pi X)\textrm{sin}(2\pi Y) \\
         Z
    \end{array}
    \right], &&\\
\end{align*}
from which we derive the values of the forcing terms $f_c$ and $\mathbf{f_u}$ and the values of the Dirichlet boundary conditions, while $\Gamma_N^c = \emptyset$. Due to this choice, the three equations are naturally decoupled ($\beta=1$ for any $t\geq0$). For this test case, we will neglect the hypothesis that the initial condition of $g = 0$. 
\par
\begin{figure}[t]
    \begin{subfigure}[b]{0.5\textwidth}
    \centering
    \resizebox{\textwidth}{!}{\definecolor{mycolor2}{rgb}{0.00000,1.00000,1.00000}%
\pgfplotsset{
  log x ticks with fixed point/.style={
      xticklabel={
        \pgfkeys{/pgf/fpu=true}
        \pgfmathparse{exp(\tick)}%
        \pgfmathprintnumber[fixed  zerofill, precision=2]{\pgfmathresult}
        \pgfkeys{/pgf/fpu=false}
      }
  }
}
\begin{tikzpicture}

\begin{axis}[%
width=3.875in,
height=2.36in,
at={(2.6in,1.099in)},
scale only axis,
xmode=log,
xmin=0.1083,
xmax=0.8660,
xminorticks=true,
xlabel = {$h$ [-]},
ylabel = {},
ymode=log,
ymin=1e-8,
ymax=5,
yminorticks=true,
axis background/.style={fill=white},
title style={font=\bfseries},
title={},
xmajorgrids,
xminorgrids,
ymajorgrids,
yminorgrids,
legend pos = south east,
legend style={legend cell align=left, align=left, draw=white!15!black}
]
              
\addplot [color=blue, line width=2.0pt]
  table[row sep=crcr]{%
0.1083  0.001288\\
0.2165  0.006196\\
0.4330  0.028615\\
0.8660  0.102231\\
};
\addlegendentry{$p=1$}

\addplot [color=mycolor2, line width=2.0pt]
  table[row sep=crcr]{%
0.1083  0.000008\\
0.2165  0.000106\\
0.4330  0.001281\\
0.8660  0.013033\\
};
\addlegendentry{$p=2$}

\addplot [color=green, line width=2.0pt]
  table[row sep=crcr]{%
0.1083  3.5791934774699415e-07\\
0.2165  0.000006\\
0.4330  0.000098\\
0.8660  0.001509\\
};
\addlegendentry{$p=3$}

\node[right, align=left, text=black, font=\footnotesize]
at (axis cs:0.159,1.3e-3) {$2$};

\addplot [color=black, line width=1.5pt]
  table[row sep=crcr]{%
0.1259   1e-3\\
0.1585   0.0016\\
0.1585   1e-3\\
0.1259   1e-3\\
};

\node[right, align=left, text=black, font=\footnotesize]
at (axis cs:0.159,1.3e-5) {$3$};

\addplot [color=black, line width=1.5pt]
  table[row sep=crcr]{%
0.1259   9e-6\\
0.1585   1.7957e-05\\
0.1585   9e-6\\
0.1259   9e-6\\
};

\node[right, align=left, text=black, font=\footnotesize]
at (axis cs:0.159,4.5e-7) {$4$};

\addplot [color=black, line width=1.5pt]
  table[row sep=crcr]{%
0.1259   3e-7\\
0.1585   7.5357e-07\\
0.1585   3e-7\\
0.1259   3e-7\\
};

\end{axis}
\end{tikzpicture}%
}
        \caption[]%
         {{\small $\|c(T)-c_h^N\|_{L^2}$}}    
        \label{fig:test2_cL2}
    \end{subfigure}\hfill
    \begin{subfigure}[b]{0.5\textwidth}
      \centering
    \resizebox{\textwidth}{!}{\definecolor{mycolor2}{rgb}{0.00000,1.00000,1.00000}%
\pgfplotsset{
  log x ticks with fixed point/.style={
      xticklabel={
        \pgfkeys{/pgf/fpu=true}
        \pgfmathparse{exp(\tick)}%
        \pgfmathprintnumber[fixed  zerofill, precision=2]{\pgfmathresult}
        \pgfkeys{/pgf/fpu=false}
      }
  }
}
\begin{tikzpicture}

\begin{axis}[%
width=3.875in,
height=2.36in,
at={(2.6in,1.099in)},
scale only axis,
xmode=log,
xmin=0.1083,
xmax=0.8660,
xminorticks=true,
xlabel = {$h$ [-]},
ylabel = {},
ymode=log,
ymin=1e-7,
ymax=50,
yminorticks=true,
axis background/.style={fill=white},
title style={font=\bfseries},
title={},
xmajorgrids,
xminorgrids,
ymajorgrids,
yminorgrids,
legend pos = south east,
legend style={legend cell align=left, align=left, draw=white!15!black}
]
              
\addplot [color=red, line width=2.0pt]
  table[row sep=crcr]{%
0.1083   0.0864186092872842\\
0.2165   0.238596\\
0.4330   0.596864\\
0.8660   1.368857 \\
};
\addlegendentry{$p=1$}

\addplot [color=orange, line width=2.0pt]
  table[row sep=crcr]{%
0.1083   0.003102\\
0.2165   0.017735\\
0.4330   0.098062\\
0.8660   0.498766 \\
};
\addlegendentry{$p=2$}

\addplot [color=purple, line width=2.0pt]
  table[row sep=crcr]{%
0.1083   9.94964241403641e-05\\
0.2165   0.001065\\
0.4330   0.011171\\
0.8660   0.113980 \\
};
\addlegendentry{$p=3$}

\node[right, align=left, text=black, font=\footnotesize]
at (axis cs:0.159,0.075) {$1$};

\addplot [color=black, line width=1.5pt]
  table[row sep=crcr]{%
0.1259   7e-2\\
0.1585   0.0881\\
0.1585   7e-2\\
0.1259   7e-2\\
};

\node[right, align=left, text=black, font=\footnotesize]
at (axis cs:0.159,3.5e-3) {$2$};

\addplot [color=black, line width=1.5pt]
  table[row sep=crcr]{%
0.1259   3e-3\\
0.1585   0.0048\\
0.1585   3e-3\\
0.1259   3e-3\\
};

\node[right, align=left, text=black, font=\footnotesize]
at (axis cs:0.159,1.4e-4) {$3$};

\addplot [color=black, line width=1.5pt]
  table[row sep=crcr]{%
0.1259   1e-4\\
0.1585   1.9953e-04\\
0.1585   1e-4\\
0.1259   1e-4\\
};

\end{axis}
\end{tikzpicture}%
}
        \caption[]%
        {{\small $\|c(T)-c_h^N\|_{DG}$}}    
        \label{fig:test2cdg}
    \end{subfigure}\hfill
    \begin{subfigure}[b]{0.5\textwidth}
      \centering
    \resizebox{\textwidth}{!}{\definecolor{mycolor2}{rgb}{0.00000,1.00000,1.00000}%
\pgfplotsset{
  log x ticks with fixed point/.style={
      xticklabel={
        \pgfkeys{/pgf/fpu=true}
        \pgfmathparse{exp(\tick)}%
        \pgfmathprintnumber[fixed  zerofill, precision=2]{\pgfmathresult}
        \pgfkeys{/pgf/fpu=false}
      }
  }
}
\begin{tikzpicture}

\begin{axis}[%
width=3.875in,
height=2.36in,
at={(2.6in,1.099in)},
scale only axis,
xmode=log,
xmin=0.1083,
xmax=0.8660,
xminorticks=true,
xlabel = {$h$ [-]},
ylabel = {},
ymode=log,
ymin=1e-8,
ymax=5,
yminorticks=true,
axis background/.style={fill=white},
title style={font=\bfseries},
title={},
xmajorgrids,
xminorgrids,
ymajorgrids,
yminorgrids,
legend pos = south east,
legend style={legend cell align=left, align=left, draw=white!15!black}
]
              
\addplot [color=blue, line width=2.0pt]
  table[row sep=crcr]{%
0.1083  0.0112760713628757\\
0.2165  0.045810\\
0.4330  0.168614\\
0.8660  0.335696\\
};
\addlegendentry{$p=1$}

\addplot [color=mycolor2, line width=2.0pt]
  table[row sep=crcr]{%
0.1083  0.000231306984591788\\
0.2165  0.002086\\
0.4330  0.020520\\
0.8660  0.133198\\
};
\addlegendentry{$p=2$}

\addplot [color=green, line width=2.0pt]
  table[row sep=crcr]{%
0.1083  2.06782386990378e-05\\
0.2165  0.000332\\
0.4330  0.004719\\
0.8660  0.079140\\
};
\addlegendentry{$p=3$}

\node[right, align=left, text=black, font=\footnotesize]
at (axis cs:0.159,1.6e-2) {$2$};

\addplot [color=black, line width=1.5pt]
  table[row sep=crcr]{%
0.1259   1.2e-2\\
0.1585   0.0190\\
0.1585   1.2e-2\\
0.1259   1.2e-2\\
};

\node[right, align=left, text=black, font=\footnotesize]
at (axis cs:0.159,4.2e-4) {$3$};

\addplot [color=black, line width=1.5pt]
  table[row sep=crcr]{%
0.1259   3e-4\\
0.1585   5.9858e-04\\
0.1585   3e-4\\
0.1259   3e-4\\
};

\node[right, align=left, text=black, font=\footnotesize]
at (axis cs:0.159,7.5e-5) {$4$};

\addplot [color=black, line width=1.5pt]
  table[row sep=crcr]{%
0.1259   3e-5\\
0.1585   7.5357e-5\\
0.1585   3e-5\\
0.1259   3e-5\\
};

\node[right, align=left, text=black, font=\footnotesize]
at (axis cs:0.1005,9e-6) {$4$};

\addplot [color=black, line width=1.5pt]
  table[row sep=crcr]{%
0.100   2.000e-05\\
0.075   6.328e-06\\
0.100   6.328e-06\\
0.100   2.000e-05\\
};

\node[right, align=left, text=black, font=\footnotesize]
at (axis cs:0.1005,4e-7) {$5$};

\addplot [color=black, line width=1.5pt]
  table[row sep=crcr]{%
0.100   8.000e-07\\
0.075   1.898e-07\\
0.100   1.898e-07\\
0.100   8.000e-07\\
};

\node[right, align=left, text=black, font=\footnotesize]
at (axis cs:0.1005,2e-8) {$6$};

\addplot [color=black, line width=1.5pt]
  table[row sep=crcr]{%
0.100   4.0000e-8\\
0.075   7.7119e-9\\
0.100   7.7119e-9\\
0.100   4.0000e-8\\
};

\end{axis}
\end{tikzpicture}%
}
        \caption[]%
        {{\small $\|\mathbf{u}(T)-\mathbf{u_h}^N\|_{L^2}$}}  
        \label{fig:test2ul2}
    \end{subfigure}\hfill
    \begin{subfigure}[b]{0.5\textwidth}
      \centering
    \resizebox{\textwidth}{!}{\definecolor{mycolor2}{rgb}{0.00000,1.00000,1.00000}%
\pgfplotsset{
  log x ticks with fixed point/.style={
      xticklabel={
        \pgfkeys{/pgf/fpu=true}
        \pgfmathparse{exp(\tick)}%
        \pgfmathprintnumber[fixed  zerofill, precision=2]{\pgfmathresult}
        \pgfkeys{/pgf/fpu=false}
      }
  }
}
\begin{tikzpicture}

\begin{axis}[%
width=3.875in,
height=2.36in,
at={(2.6in,1.099in)},
scale only axis,
xmode=log,
xmin=0.1083,
xmax=0.8660,
xminorticks=true,
xlabel = {$h$ [-]},
ylabel = {},
ymode=log,
ymin=1e-5,
ymax=3000,
yminorticks=true,
axis background/.style={fill=white},
title style={font=\bfseries},
title={},
xmajorgrids,
xminorgrids,
ymajorgrids,
yminorgrids,
legend pos = south east,
legend style={legend cell align=left, align=left, draw=white!15!black}
]
              
\addplot [color=red, line width=2.0pt]
  table[row sep=crcr]{%
0.1083   30.9731554523587\\
0.2165   66.656590\\
0.4330   142.480244\\
0.8660   204.713471\\
};
\addlegendentry{$p=1$}

\addplot [color=orange, line width=2.0pt]
  table[row sep=crcr]{%
0.1083   1.95094374828515\\
0.2165   9.458818\\
0.4330   44.860178\\
0.8660   156.840452\\
};
\addlegendentry{$p=2$}

\addplot [color=purple, line width=2.0pt]
  table[row sep=crcr]{%
0.1083   0.1826174896770226\\
0.2165   1.420003\\
0.4330   11.674007\\
0.8660   105.375572\\
};
\addlegendentry{$p=3$}

\node[right, align=left, text=black, font=\footnotesize]
at (axis cs:0.159,33) {$1$};

\addplot [color=black, line width=1.5pt]
  table[row sep=crcr]{%
0.1259   3e1\\
0.1585   37.7678\\
0.1585   3e1\\
0.1259   3e1\\
};

\node[right, align=left, text=black, font=\footnotesize]
at (axis cs:0.159,2.5) {$2$};

\addplot [color=black, line width=1.5pt]
  table[row sep=crcr]{%
0.1259   2\\
0.1585   3.1698\\
0.1585   2\\
0.1259   2\\
};

\node[right, align=left, text=black, font=\footnotesize]
at (axis cs:0.159,0.28) {$3$};

\addplot [color=black, line width=1.5pt]
  table[row sep=crcr]{%
0.1259   0.2\\
0.1585   0.3991\\
0.1585   0.2\\
0.1259   0.2\\
};

\end{axis}
\end{tikzpicture}%
}
        \caption[]%
        {{\small $\|\mathbf{u}(T)-\mathbf{u_h}^N\|_{DG}$}}   
        \label{fig:test2udg}
\end{subfigure}
\caption{Test case 1: computed $L^2$ (left) and $DG$ (right) errors as functions of the mesh-size, with $h = 0.8660, 0.4330, 0.2165, 0.1083$ and $p = 1, 2, 3$, for the test case 1. The top panels show the results for the approximate concentration, while bottom panels show the results for the approximate displacement. All the errors are computed at final time $T = 0.1$ years. Logarithmic scales are used for all axes. The triangles indicate different convergence rates.}
\label{fig:test2_errors}
    \end{figure}
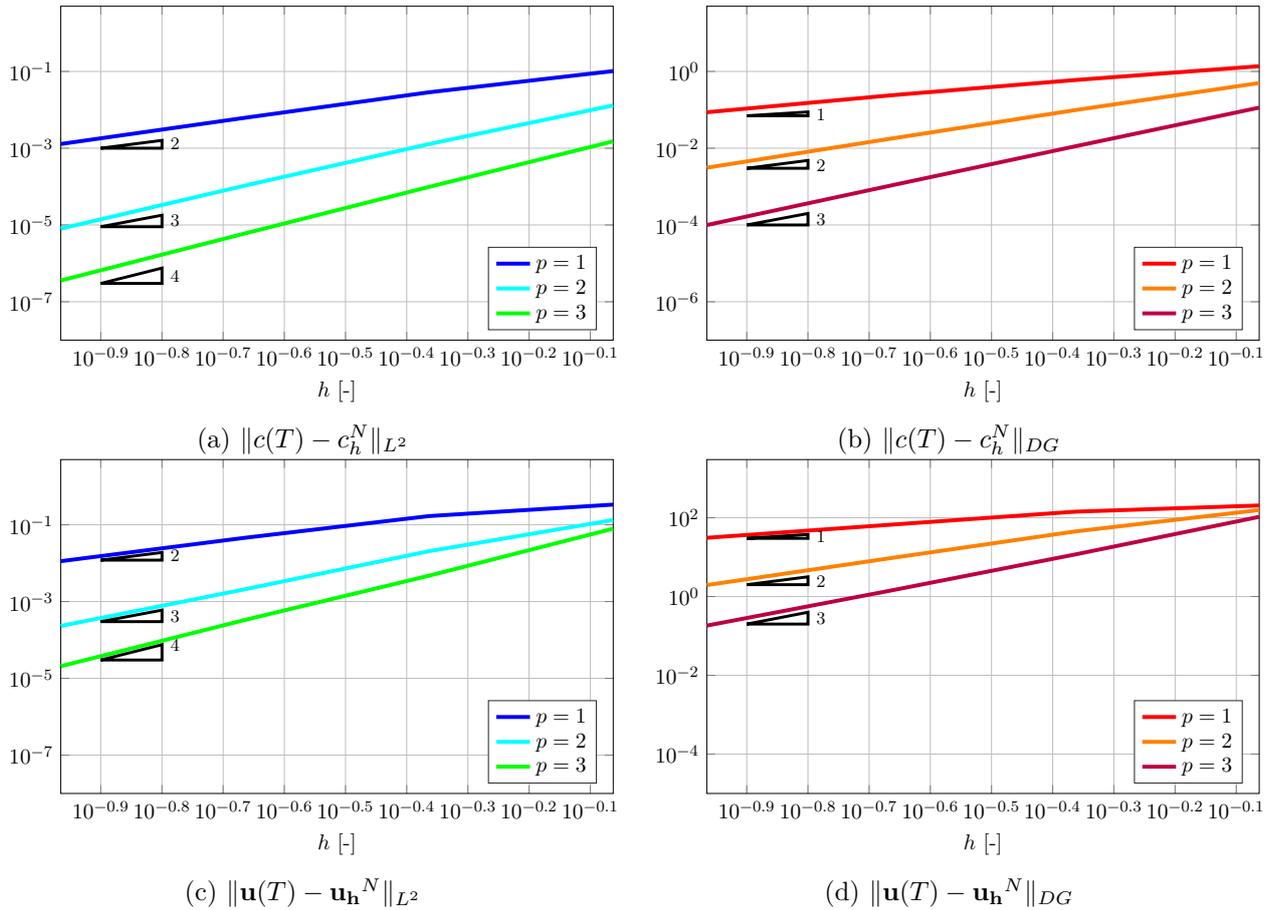

Figure \ref{fig:test2_errors} shows the computed errors in the $L^2$ and $DG$ norms of the errors as functions of the mesh size $h = 0.8660, \ 0.4330, \ 0.2165, \ 0.1083$, in logarithmic scale. We can see how the error’s norms follow the expected trend of $h^{p+1}$ for the $L^2$-norms and of $h^p$ for the $DG$-norms as proven for the FK equation, and the linear elasticity equation \cite{FK_equation_and_PolyDG, Poroelasticity_PolyDG}.

\subsection{Test Case 2: coupled system}
\label{testcase2}
In this section, we present the results obtained for a holed spherical domain $\Omega = \{\mathbf{X} \in \mathbb{R}^3: 0.05^2 \le X^2 + Y^2 + Z^2 \le 0.1^2\}$. We define the outer boundary as the Neumann's boundary $\Gamma_N:= \{\mathbf{X}\in \mathbb{R}^3:\,X^2 + Y^2 + Z^2 = 0.1^2 \}$ and the inner boundary as the Dirichlet's boundary: $\Gamma_D:= \{\mathbf{X}\in \mathbb{R}^3:\,X^2 + Y^2 + Z^2 = 0.05^2 \}$.
\par
\begin{figure}[t]
    \centering
         \resizebox{\textwidth}{!}
        {\includegraphics[scale=0.08]{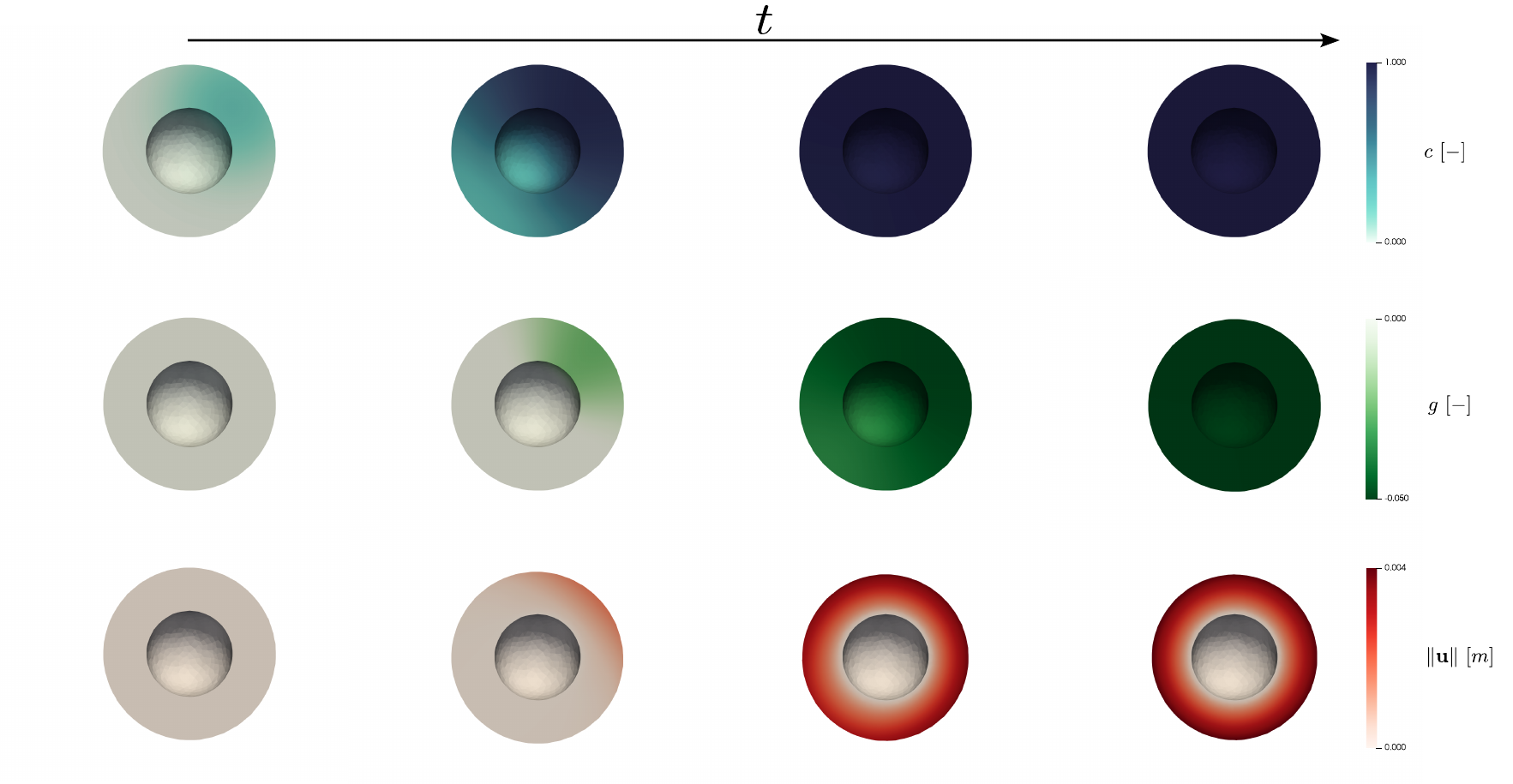}}
    \caption
    {Test Case 2: concentration, atrophy rate, and norm of the displacement with $\widehat{\tau} = 1$ year at time $t = 0, 5, 10, 15$ years, from left to right, obtained with polynomial degree $p = 2$.} 
    \label{fig:1simulation}
\end{figure}

Again, we employ the parameters presented in Table \ref{tab:test2_parameters} and consider $f_c = 0$. As for the initial condition, $c_0(\mathbf{X})$, we assume the pathogen is primarily concentrated in a small portion of the domain. Specifically, we use a Gaussian-type distribution to represent this initial state
\begin{equation}
    c(\mathbf{X}, 0) = A\exp{\left[-\frac{1}{2}\left(\frac{|\mathbf{X} - \mathbf{X_0}|}{0.05}\right)^2\right]},
\end{equation}
where $A \simeq 0.4$ represents the amplitude, and the center $\mathbf{X_0}$ is chosen to be $(0.05,0.05,0.05)$.
\begin{figure}
    \centering
         \resizebox{\textwidth}{!}
        {\includegraphics[scale=0.08]{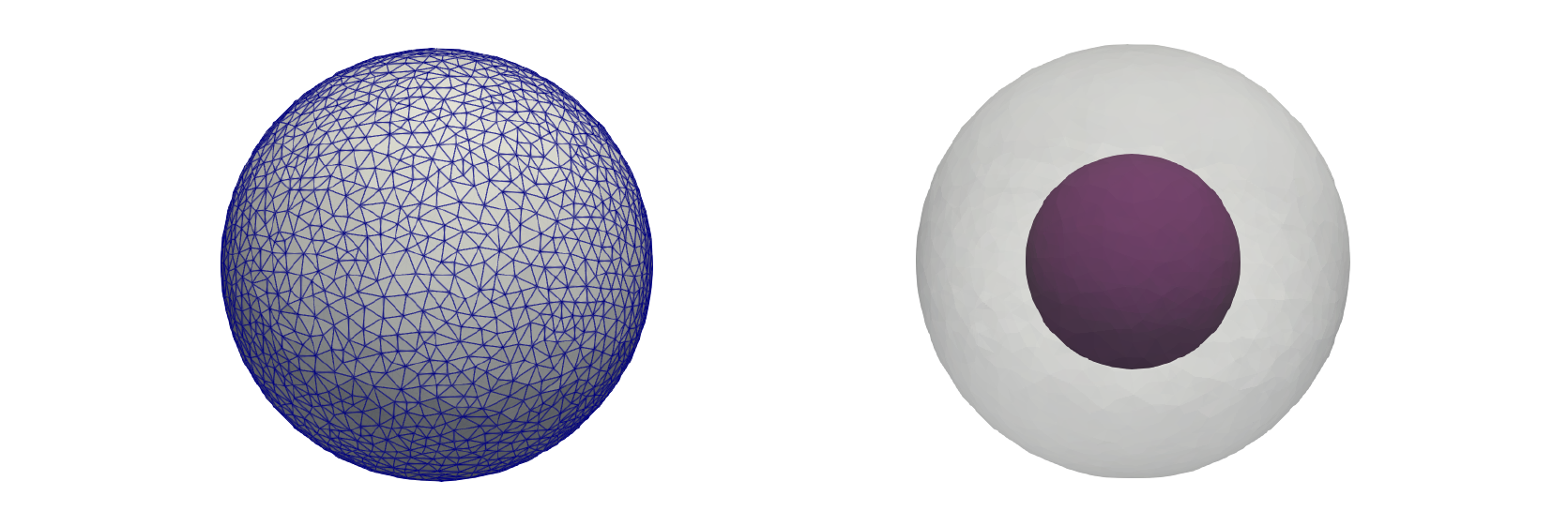}}
    \caption
    {An external representation of the mesh (left) and a visualization of the ventricles boundary in purple and of the skull in transparency (right).} 
    \label{fig:mesh sfera}
\end{figure}.
Finally we enforce homogeneous Neumann boundary conditions, to avoid the spreading of the pathogen outside the domain. For the displacement we impose we impose homogeneous Neumann boundary conditions on $\Gamma_N$ (the outer boundary), as we do not have any external stress, and homogeneous Dirichlet boundary conditions on $\Gamma_D$ (inner boundary), to block the infinitesimal rigid displacement. 
The solutions in Figure \ref{fig:1simulation} are computed with mesh size $h = 0.094$, time interval $\Delta t=0.05$ years, and polynomial degree $p=3$. We show a section of the holed sphere where we can see the initial concentration distribution and its evolution.  At final time $T = 15 \ [\text{years}]$, we can see how both the concentration and the length reduction rate $g$ have reached their maximal value, $c_h(\mathbf{X}, T) = 1.0$ and, $|g_h(\mathbf{X}, T)| = \gamma = 0.05$, in agreement with what is expected from the model.

\section{Simulations on a real brain geometry}\label{sec: Simulations on brain}
In this section, we present the simulation by applying our model to the problem of the onset of brain atrophy induced by Alzheimer's disease.
In the following, we assume that the brain undergoes only small deformations. It is important to note that the brain is composed of extremely soft material, even compared to other biological soft tissues \cite{budday2017mechanical}. Although nonlinear elasticity might initially appear to be a more appropriate framework for modeling such soft tissues, in this case, the mechanical deformations primarily represent the atrophy induced by Alzheimer's disease. Throughout the disease, the brain's atrophy rate is approximately $1\%$-$2\%$ per year for the entire brain \cite{Sch_fer_2021}. Consequently, the resulting atrophy remains relatively small, even after a decade, with an estimated mean value of $g$ between $-0.032$ and $-0.065$. We perform our simulations on the real geometry of a brain derived from the MRI images of project OASIS \cite{OASIS}. For this simulation, we consider a transversely isotropic diffusion tensor in the form $\mathbf{D} = d_{ext}\mathbf{I} + d_{axn} (\mathbf{n}\otimes\mathbf{n})$, and the parameters listed in table \eqref{tab:testreal_parameters}.

\begin{table}[t]
\caption{Test Case 3: Parameters used in the Test Case 3}\label{tab:testreal_parameters}
\begin{tabular}{lrlc|lrlc}
\toprule
  \textbf{Parameters} & \multicolumn{2}{c}{\textbf{Values}} &
  \textbf{Reference} & \textbf{Parameters} & \multicolumn{2}{c}{\textbf{Values}} &
  \textbf{Reference} \\ 
\midrule
 \normalsize$d_{ext}$ & \normalsize$8.00$ & $[\mathrm{mm^2/year}]$ &
 \normalsize\cite{biochemical_biomechanical_Kuhl} & \normalsize$\lambda$ &
 \normalsize $2700$ & $[\mathrm{Pa}]$&
 \normalsize \cite{brainparam1, brainparam2}\\
 \normalsize$d_{axn}$ & \normalsize $80.00 $ & $\mathrm{[mm^2/year]}$ &
 \normalsize \cite{biochemical_biomechanical_Kuhl} & \normalsize $\mu$ & \normalsize $300$ & $[\mathrm{Pa}]$ &
 \normalsize \cite{brainparam1, brainparam2}\\
 \normalsize $\alpha$ &
 \normalsize $0.90$ & $[-]$&
 \normalsize \cite{biochemical_biomechanical_Kuhl} & \normalsize $\widehat{\tau}$ & \normalsize $1$ & $[\mathrm{year}]$ &
 \normalsize  \\
  \normalsize $\gamma$ & \normalsize $0.05$ & $[-]$ &
 \normalsize  \\
\bottomrule
\end{tabular}
\end{table}
\par
\begin{figure}
    \centering
         \resizebox{\textwidth}{!}
        {\includegraphics[scale=0.08]{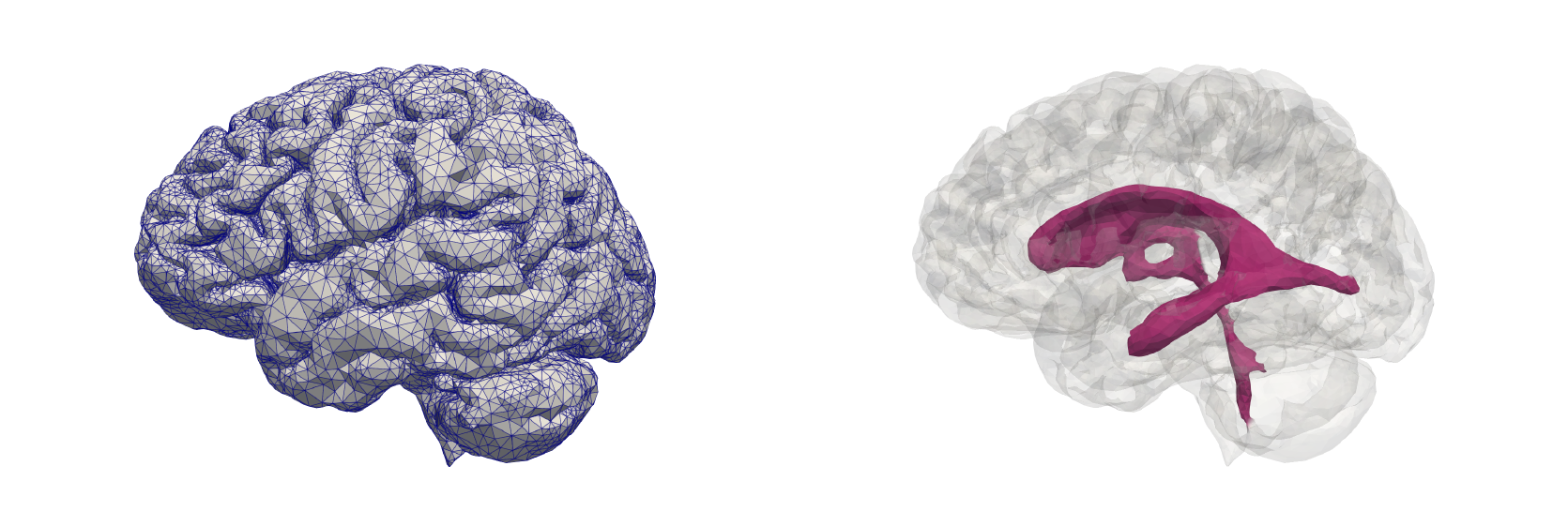}}
    \caption
    {An external representation of the mesh (left) and a visualization of the ventricles boundary in purple and of the skull in transparency (right).} 
    \label{fig:mesh cervello}
\end{figure}

We impose homogeneous Neumann boundary conditions for the concentration on the whole boundary $\partial\Omega$. Concerning the displacement, we impose homogeneous Dirichlet conditions on the ventricular boundary $\Gamma_{\text{Vent}}$, and homogeneous Neumann conditions on the skull $\Gamma_{\text{Skull}}$ (see Figure \ref{fig:mesh cervello}).
 The discretization is performed with polynomial degree $p=2$ for the concentration $c$, and polynomial degree $p=3$ for the atrophy rate $g$ and the displacement $\mathbf{u}$. Additionally, we apply a temporal discretization with $\Delta t = 0.05$ years.
\begin{figure}
    \centering
         \resizebox{\textwidth}{!}
        {\includegraphics[scale=0.08]{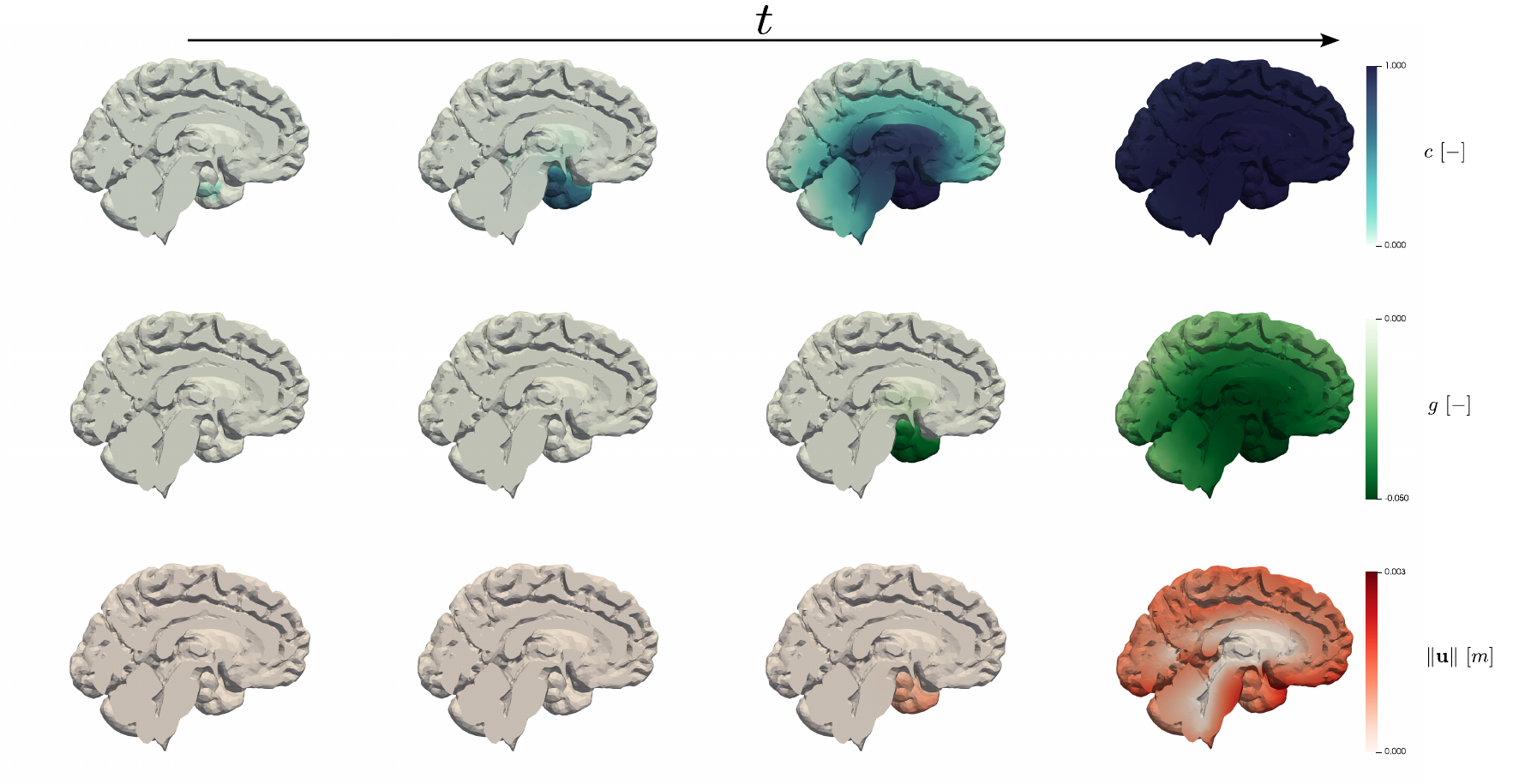}}
    \caption
    {Test Case Section \ref{sec: Simulations on brain} Concentration, atrophy rate and displacement with $\widehat{\tau} = 1$ year at time $t = 0, 6, 12, 18$ years obtained with polynomial degree 2 for concentration and 3 for atrophy rate and displacement.} 
    \label{fig:2simulation}
\end{figure}
\par
The results of the simulations are reported in Figure \ref{fig:2simulation}, which shows us the results for times $t = 0, 6, 12, 18$ years in a section of the brain. We can observe the presence of a seeding of misfolded $\tau$-proteins located in the locus coeruleus \cite{juckerSelfpropagationPathogenicProtein2013}. The proteins spread in all directions according to the axonal directions. The variable $g$ activates later, and the displacement field varies accordingly. At the final time, $t = 18$, the concentration of proteins reaches its maximum value in the whole domain, and we can observe a shrinkage of the cortical surface.

\section{Extension with nonlinear elasticity}\label{sec: extension nonlinear}

\begin{table}
\caption{Test Case 4: Parameters used in the test case described in sections~\ref{sec: extension nonlinear}}\label{tab:testreal_parameters_cube}
\begin{tabular}{lrlc|lrlc}
\toprule
  \textbf{Parameters} & \multicolumn{2}{c}{\textbf{Values}} &
  \textbf{Reference} & \textbf{Parameters} & \multicolumn{2}{c}{\textbf{Values}} &
  \textbf{Reference} \\ 
\midrule
 \normalsize$d_{ext}$ & \normalsize$8.00$ & $[\mathrm{mm^2/years}]$ &
 \normalsize\cite{biochemical_biomechanical_Kuhl} & \normalsize$\lambda$ &
 \normalsize $505$ & $[\mathrm{Pa}]$&
 \normalsize \cite{Poroelasticity_PolyDG}\\
 \normalsize$d_{axn}$ & \normalsize $80.00 $ & $\mathrm{[mm^2/years]}$ &
 \normalsize \cite{biochemical_biomechanical_Kuhl} & \normalsize $\mu$ & \normalsize $216$ & $[\mathrm{Pa}]$ &
 \normalsize \cite{Poroelasticity_PolyDG}\\
 \normalsize $\alpha$ &
 \normalsize $0.90$ & $[-]$&
 \normalsize \cite{biochemical_biomechanical_Kuhl} & \normalsize $\widehat{\tau}$ & \normalsize $1$ & $[\mathrm{year}]$ &
 \normalsize  \\
  \normalsize $\gamma$ & \normalsize $0.05$ & $[-]$ &
 \normalsize  \\
\bottomrule
\end{tabular}
\end{table}
We finally extend our model in the more challenging setting of nonlinear elasticity. We consider the following hyperelastic energy:
\begin{equation}
    \Psi_0(\mathbf{F}_E) = \frac{\mu}{2} (\mathbf{F}_E:\mathbf{F}_E - 3 - 2\log(\det\mathbf{F}_E)) + \frac{\lambda}{2}(\det\mathbf{F}_E - 1)^2.
\end{equation}
As before, we consider the tensor $\mathbf{G} = (1+g)\mathbf{I}$. Since small deformations are no longer necessary, we can set the parameter $\gamma$ (representing the maximum volume loss when the concentration of misfolded tau proteins equals 1) to reach values around $20\%$. In this case, we have a variation of the volume of the order of $O((1+g)^3)$. Under such hypothesis we cannot assume that the current configuration $\Omega_t \equiv \Omega$, the reference configuration.
\par
\begin{figure}
    \centering
         \resizebox{0.8\textwidth}{!}
        {\includegraphics[scale=0.08]{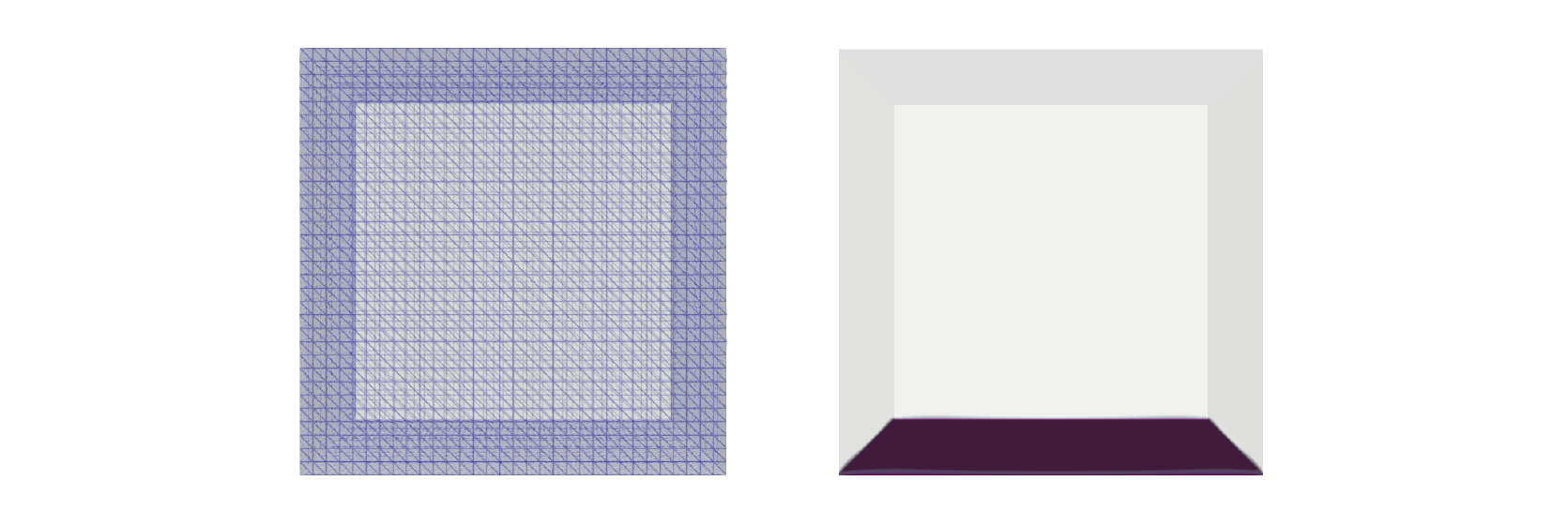}}
    \caption
    {An external representation of the mesh (left) and a visualization of the Dirichlet's boundary in purple and of the Neumann's boundary in transparency (right).} 
    \label{fig:meshBCcubo}
\end{figure}

We perform the simulations using discontinuous elements for the FK and the logistic equation and continuous elements for the nonlinear elasticity equation. In this section, we present the results of numerical simulations on a three-dimensional tetrahedral grid. We consider a cubic domain, denoted as $\Omega = (0\,\mathrm{m },1\,\mathrm{m})^3$, and the parameters detailed in \cite{FK_equation_and_PolyDG, Poroelasticity_PolyDG} for the concentration and elasticity equations, as presented in Table \ref{tab:testreal_parameters_cube}. Additionally, we enforce a $c_{cr} = 0.8$. Regarding the concentration, we introduce a zero forcing term $f_c = 0$ and homogeneous Neumann boundary conditions.
Moreover, we define the initial condition as a Gaussian function with an amplitude $A$ of approximately 0.8, ensuring the concentration starts below the critical value:
\begin{equation}
    c(\mathbf{X}, 0) = A\exp{\left[-\frac{1}{2}\left(\frac{|\mathbf{X}-\mathbf{X_0}|}{0.15}\right)^2\right]}
\end{equation}
We take $g(\mathbf{X},\,0)=0$ as the initial condition for the atrophy field. We further assume that the body is free of external body forces (i.e., $\mathbf{f_u} = \mathbf{0}$). Finally, we impose homogeneous Dirichlet’s boundary conditions on $\Gamma^\mathbf{u}_D = [0, 1] \times [0, 1] \times \{Z = 1\}$ along with homogeneous Neumann’s conditions on $\Gamma^\mathbf{u}_N =\partial\Omega\setminus\Gamma^\mathbf{u}_D$, Figure \ref{fig:meshBCcubo}.  
\begin{figure}
    \centering
         \resizebox{\textwidth}{!}
        {\includegraphics[scale=0.08]{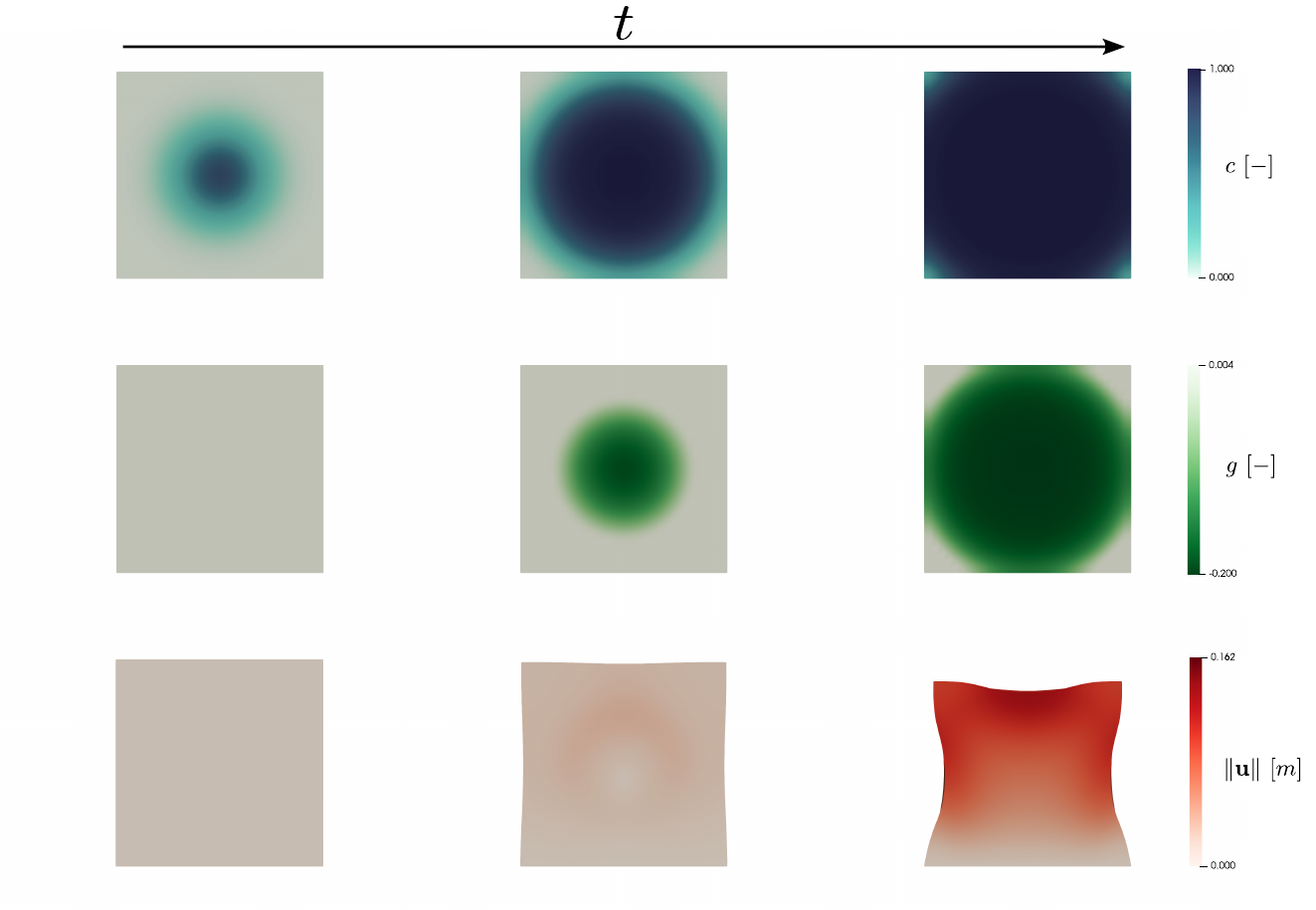}}
    \caption
    {Test Case Section \ref{sec: extension nonlinear}. Concentration, atrophy rate and displacement with $\widehat{\tau} = 1$ year at time $t = 0, 5, 10$ years obtained with polynomial degree 2 for concentration and 3 for atrophy rate and displacement.} 
    \label{fig:3simulation}
\end{figure}
The results are shown in Figure \ref{fig:3simulation}. As expected, we observe that the cube shrinks accordingly, with the misfolded proteins spread through the atrophy rate.

\section{Conclusion} \label{sec:conclusion}
In this work, we have presented a multiphysics model for studying tissue atrophy caused by pathogen spread, with a specific application to Alzheimer's disease. The pathogen's spread and aggregation have been modeled using the Fisher-Kolmogorov equation, incorporating both dispersion and proliferation effects. Tissue atrophy has been described through a morpho-elastic framework, where mass loss and tissue elasticity together shape the resulting tissue morphology. The model integrates the morpho-elastic response with pathogen concentration by introducing an evolution law for inelastic strain, governed by pathogen concentration through a logistic-type differential equation.

For the construction of the discrete model, we have employed a PolyDG method, exploiting its versatility in the discretization of complex domains. We considered a Crank-Nicolson method for the time derivatives discretization and a semi-implicit treatment of the non-linear term. We have assessed the validity of the PolyDG method via two convergence tests, discussing the results with respect to the theoretical outcomes. Secondly, we have validated our model by conducting simulations on a holed-spherical domain. Moreover, we have applied the model to the simulation of Alzheimer's disease on a real brain geometry, where we have observed outcomes consistent with the anticipated biological behavior of prion-like protein diffusion and tissue atrophy. While the small elasticity assumption for the whole brain holds due to the relatively low mean value of the atrophy parameter $g$, some brain regions may experience more significant volume loss. As a proof-of-concept, we also performed a simulation incorporating a nonlinear constitutive law for tissue elasticity. In this case, we used discontinuous elements for both the pathogen concentration and the atrophy variable $g$, while opting for continuous elements to discretize the displacement field for simplicity. Possible future works include the implementation of more precise models of the diffusion of the pathogen, such as the heterodimer model.

\subsection*{Acknowledgments}
The brain MRI images were provided by OASIS-3: Longitudinal Multimodal Neuroimaging: Principal Investigators: T. Benzinger, D. Marcus, J. Morris; NIH P30 AG066444, P50 AG00561, P30 NS09857781, P01 AG026276, P01 AG003991, R01 AG043434, UL1 TR000448, R01 EB009352. AV-45 doses were provided by
Avid Radiopharmaceuticals, a wholly-owned subsidiary of Eli Lilly.

\section*{Appendix A: FK equation in reference configuration}\label{AppB}
In the current configuration $\Omega$ the Fisher-Kolmogorov equation for the concentration of misfolded proteins reads:
\begin{equation}
    \frac{\partial c}{\partial t} = \nabla \cdot (\mathbf{D}\nabla c) + \alpha c(1-c) + f_c.
\end{equation}
The integral of a function $f$ in the current configuration is equal to the integral of the same function multiplied for the determinant of the deformation gradient in the reference configuration: 
\begin{equation}\label{eq: forref}
    \int_{\mathcal{P}_t} f d\mathbf{x} = \int_{\mathcal{P}_0} J f d\mathbf{X}.
\end{equation}
Now we consider the partial derivatives:
\begin{equation}
    \frac{\partial}{\partial x_i} = \frac{\partial}{\partial X_j} \frac{\partial X_j}{\partial x_i}
\end{equation}
We define $(\nabla_X)_j$ as $\frac{\partial}{\partial X_j}$, while $ \frac{\partial X_j}{\partial x_i}$ is the $ij$ component of the inverse of the deformation gradient $F^{-1}$. Indeed, we obtain that:
\begin{equation}\label{eq: gradref}
    \nabla = \mathbf{F}^{-T} \nabla_X.
\end{equation}
Now we need to define the divergence operator ($\nabla \cdot$ ) in the reference configuration.
We consider the volume integral over a portion $\mathcal{P_t}$ of the current configuration of the divergence of a function $\mathbf{v}$ and we apply the divergence theorem:
\begin{equation}
    \int_{\mathcal{P}_t} \nabla \cdot \mathbf{v} \ d\mathbf{x} = \int_{\partial \mathcal{P}_t} \mathbf{v}\cdot\mathbf{n} \ ds.
\end{equation}
Now we apply the  Nanson's formula:
\begin{equation}
    \int_{\partial \mathcal{P}_t} \mathbf{v}\cdot\mathbf{n} \ ds = \int_{\partial \mathcal{P}_0} \mathbf{v}\cdot (J \mathbf{F}^{-T}\mathbf{N}) \ dS,
\end{equation}
then we reapply the divergence theorem:
\begin{equation} \label{eq: divref}
    \int_{\partial \mathcal{P}_0} \mathbf{v}\cdot (J \mathbf{F}^{-T}\mathbf{N}) \ dS = \int_{\partial \mathcal{P}_0}  (J \mathbf{F}^{-1}\mathbf{v}) \cdot \mathbf{N} \ dS =  \int_{\mathcal{P}_0} \nabla_X \cdot (J \mathbf{F}^{-1}\mathbf{v}) \ dS.
\end{equation}
Now applying  \eqref{eq: forref}, \eqref{eq: gradref}, \eqref{eq: divref} to the FK equation we obtain 
\begin{equation}
    J\frac{\partial c}{\partial t} = \nabla_X \cdot (J\mathbf{F}^{-1}\mathbf{D}\mathbf{F}^{-T}\nabla_X c) + \alpha Jc(1-c).
\end{equation}

\section*{Appendix B: Matrices' definitions}\label{AppA}
Employing the same basis defined in Section \ref{algebraicModel}, we can define the following matrices for $i,j = 1,...,N_c$:
\begin{align}
    &[\mathbf{M}_c]_{ij} =  (\phi_j, \phi_i)_\Omega & (\mathrm{Mass\;matrix}); \\ &[\mathbf{A}_c]_{ij} =  \mathscr{A}_c(\phi_j, \phi_i) & (\mathrm{Stiffness\;matrix}); \\
   &[\mathbf{M}_\alpha]_{ij} =  (\alpha \phi_j, \phi_i)_\Omega & (\mathrm{Linear\;reaction\;matrix}); \\ &[\widetilde{\mathbf{M}}_\alpha(\mathbf{C}(t))]_{ij} =  (\alpha c_h(t)\phi_j, \phi_i)_\Omega & (\mathrm{Nonlinear\;reaction\;matrix});
\end{align}
We also define for $i,j = 1,...,N_g$:
\begin{align}
    &[\mathbf{M}_g]_{ij} =  (q_j,q_i)_\Omega & (g\;\mathrm{-mass\;matrix}); \\ 
   &[\mathbf{M}_\beta]_{ij} =  r_L^g(q_j,q_i) & (\mathrm{Linear\;term\;matrix}); \\
   &[\widetilde{\mathbf{M}}_\beta(\mathbf{g}(t))]_{ij} =  r_N^g(q_j, g_h(t), q_i) & (\mathrm{Nonlinear\;term\;matrix});
\end{align}
and finally set, for $i,j = 1,...,N_\mathbf{u}$:
\begin{align}
   &[\mathbf{M}_E]_{ij} =  (\bm{\psi}_j,\bm{\psi}_i)_\Omega & (\mathrm{Elasticity\;mass\;matrix}); \\ 
   &[\mathbf{K}_E]_{ij} =  \mathscr{A}_E(\bm{\psi}_j, \bm{\psi}_i) & (\mathrm{Elasticity\;stiffness\;matrix}); \\
   &[\mathbf{B}_g]_{ij} =  \mathscr{B}_E(q_j, \bm{\psi}_i) & (g-\mathrm{displacement\;coupling\;matrix}).
\end{align}
Moreover we define the forcing terms $\mathbf{F}_c = [F_c(\phi_j)]_{j=1}^{N_c}$, $\mathbf{F}_g = [F_g(q_j)]_{j=1}^{N_g}$ and $\mathbf{F}_E = [F_E(\bm{\psi}_j)]_{j=1}^{N_\mathbf{u}}$.

\bibliographystyle{hieeetr}
\bibliography{bibliografia}

\begin{thebibliography}{10}

\bibitem{physicsbasedmodelAlz}
J.~Weickenmeir, M.~Jucker, A.~Goriely, and E.~Kuhl, ``A physics-based model explains the prion-like features of neurodegeneration in alzheimer's disease, parkinson's disease, and amyotrophic lateral schlerosis'', {\em {Journal of the Mechanics and Physics of Solids}}, vol.~124, pp.~264--281, 2019.

\bibitem{FKtumor}
J.~Belmonte-Betia, G.~Calvo, and V.~Pérez-Garcia, ``Effective particle methods for {Fisher-Kolmogorov} equations: {Theory} and applications to brain tumor dynamics'', {\em {Communication in Nonlinear Science and Numerical Simulations}}, vol.~19, no.~9, pp.~3267--3283, 2014.

\bibitem{juckerSelfpropagationPathogenicProtein2013}
M.~Jucker and L.~C. Walker, ``Self-propagation of pathogenic protein aggregates in neurodegenerative diseases'', {\em Nature}, vol.~501, no.~7465, pp.~45--51, 2013.

\bibitem{ActaNeuro}
H.~Braak and E.~Braak, ``Neuropathological stages of {Alzheimer}-related changes'', {\em Acta Neuropathologica}, vol.~82, no.~4, pp.~239--259, 1991.

\bibitem{propspread}
M.~Jucker and L.~C. Walker, ``Propagation and spread of pathogenic protein assemblies in neurodegenerative diseases'', {\em Nature Neuroscience}, vol.~21, pp.~1341--1349, 2018.

\bibitem{selfprop}
M.~Jucker and L.~C. Walker, ``Self-propagation of pathogenic protein aggregates in neurodegenerative diseases'', {\em Nature}, vol.~501, no.~7465, pp.~45--51, 2013.

\bibitem{schiesser_parkinson}
W.~Schiesser, {\em {ODE/PDE} $\alpha$-synuclein models for {P}arkinson's disease}.
\newblock Elsevier, 1~ed., 2018.

\bibitem{Prionlike}
S.~Fornari, A.~Sch\"afer, M.~Jucker, A.~Goriely, and E.~Kuhl, ``Prion-like spreading of {Alzheimer}'s disease within the brain connectome'', {\em Journal of The Royal Society Interface}, vol.~16, no.~159, p.~20190356, 2019.

\bibitem{FK_equation_and_PolyDG}
M.~Corti, F.~Bonizzoni, L.~Dedé, A.~M. Quarteroni, and P.~F. Antonietti, ``Discontinuous {Galerkin} methods for {Fisher}-{Kolmogorov} equation with $\alpha$-synuclein spreading in {Parkinson}'s disease'', {\em Computer Methods in Applied Mechanics and Engineering}, vol.~417, p.~116450, 2023.

\bibitem{Heterodimer}
P.~Antonietti, F.~Bonizzoni, M.~Corti, and A.~Dall'Olio, ``{Discontinuous Galerkin} approximations of the heterodimer model for protein-protein interaction'', {\em Computer Methods in Applied Mechanics and Engineering}, vol.~431, p.~117282, 2024.

\bibitem{smoluchowski_kuhl}
S.~Fornari, A.~Sch\"afer, E.~Kuhl, and A.~Goriely, ``Spatially-extended nucleation-aggregation-fragmentation models for the dynamic of prion-like neurodegenerative protein-spreading in the brain and its connectone'', {\em Journal of Theoretical Biology}, vol.~486, p.~110102, 2020.

\bibitem{biochemical_biomechanical_Kuhl}
A.~Sch\"afer, J.~Weickenmeier, and E.~Kuhl, ``The interplay of biochemical and biomechanical degeneration in {Alzheimer}'s' disease'', {\em Computer Methods in Applied Mechanics and Engineering}, vol.~352, pp.~369--388, 2019.

\bibitem{multiph}
J.~Weickenmeier, E.~Kuhl, and A.~Goriely, ``Multiphysics of prionlike diseases: Progression and atrophy'', {\em Physical Review Letters}, vol.~121, no.~15, p.~158101, 2018.

\bibitem{Goriely_Math_Mec}
A.~Goriely, {\em The Mathematics and Mechanics of Biological Growth}.
\newblock Springer New York, NY, 2017.

\bibitem{Ambrosimorpho}
D.~Ambrosi, G.~Ateshian, A.~E.M., S.~C. Cowin, J.~Dumais, A.~Goriely, G.~Holzapfel, J.~Humphrey, E.~Kemkemer, R.~Kuhl, J.~Olberding, L.~Taber, and K.~Garikipati, ``Perspectives on biological growth and remodeling'', {\em Journal of Mechanics and Physics of Solids}, vol.~59, no.~4, pp.~863--883, 2011.

\bibitem{polydg1}
F.~Bassi, L.~Botti, A.~Colombo, D.~A. Di~Pietro, and P.~Tesini, ``On the flexibility of agglomeration based physical space discontinuous {Galerkin} discretizations'', {\em Journal of Computational Physics}, vol.~231, pp.~45--65, 2012.

\bibitem{polydg3}
J.~S. Hesthaven and T.~Warburton, {\em Nodal Discontinuous Galerkin Method: Algorithms, Analysis and Applications}.
\newblock Springer International Publishing, 2008.

\bibitem{polydg4}
P.~F. Antonietti, C.~Facciolà, P.~Houston, I.~Mazzieri, G.~Pennesi, and M.~Verani, ``High–order {Discontinuous} {Galerkin} {Methods} on {Polyhedral} {Grids} for {Geophysical} {Applications}: {Seismic} {Wave} {Propagation} and {Fractured} {Reservoir} {Simulations}'', in {\em Polyhedral {Methods} in {Geosciences}}, pp.~159--225, Springer International Publishing, 2021.

\bibitem{polydg5}
A.~Cangiani, Z.~Dong, E.~Georgoulis, and P.~Houston, {\em {hp-Version} {D}iscontinuous {G}alerkin {M}ethods on {P}olygonal and {Polyhedral} {Meshes}}.
\newblock Springer, 2017.

\bibitem{WHO}
WHO, ``Dementia'', 2023.
\newblock \url{https://www.who.int/news-room/fact-sheets/detail/dementia}.

\bibitem{Abetarate}
V.~L. Villemagne, S.~Burnham, P.~Bourgeat, B.~Brown, K.~A. Ellis, O.~Salvado, C.~Szoeke, S.~L. Macaulay, R.~Martins, P.~Maruff, D.~Ames, C.~C. Rowe, and C.~L. Masters, ``Amyloid $\beta$ deposition, neurodegeneration, and cognitive decline in sporadic alzheimer’s disease: a prospective cohort study'', {\em The Lancet, Neurology}, vol.~12, no.~4, pp.~357--367, 2013.

\bibitem{Fisher}
R.~Fisher, ``The wave of advance of advantageous genes'', {\em Annals of Eugenics}, vol.~7, no.~4, pp.~353--369, 1937.

\bibitem{Kolmogorov}
A.~N. Kolmogorov, I.~G. Petrovskii, and N.~S. Piskunov, ``Etude de la diffusion avec croissance de la quantité de matière et son application à un problème biologique'', {\em Moscow University Mathematics Bulletin}, vol.~1, pp.~1--25, 1937.

\bibitem{PDE}
S.~Salsa, {\em Equazioni a derivate parziali}.
\newblock Springer, 2016.

\bibitem{rodriguez1994stress}
E.~K. Rodriguez, A.~Hoger, and A.~D. McCulloch, ``Stress-dependent finite growth in soft elastic tissues'', {\em Journal of biomechanics}, vol.~27, no.~4, pp.~455--467, 1994.

\bibitem{gurtin2010mechanics}
M.~E. Gurtin, E.~Fried, and L.~Anand, {\em The mechanics and thermodynamics of continua}.
\newblock Cambridge university press, 2010.

\bibitem{PartialDiffEq}
S.~Salsa and G.~Verzini, {\em Partial Differential Equations in Action: From Modelling to Theory}.
\newblock Springer, 2022.

\bibitem{NumModels}
Q.~A. M., {\em Numerical Models for Differential Problems}.
\newblock Springer, 2016.

\bibitem{AlnaesEtal2015}
M.~S. Alnaes, J.~Blechta, J.~Hake, A.~Johansson, B.~Kehlet, A.~Logg, C.~N. Richardson, J.~Ring, M.~E. Rognes, and G.~N. Wells, ``The {FEniCS} project version 1.5'', {\em Archive of Numerical Software}, vol.~3, no.~100, 2015.

\bibitem{Poroelasticity_PolyDG}
M.~Corti, L.~Dedé, A.~M. Quarteroni, and P.~F. Antonietti, ``Numerical {Modelling} of the {Brain} {Poromechanics} by {Higher}-{Order} {Discontinuous} {Galerkin} {Methods}'', {\em Mathematical Models and Methods in Applied Sciences}, vol.~33, no.~8, pp.~1577--1609, 2023.

\bibitem{budday2017mechanical}
S.~Budday, G.~Sommer, C.~Birkl, C.~Langkammer, J.~Haybaeck, J.~Kohnert, M.~Bauer, F.~Paulsen, P.~Steinmann, E.~Kuhl, and G.~A. Holzapfel, ``Mechanical characterization of human brain tissue'', {\em Acta biomaterialia}, vol.~48, pp.~319--340, 2017.

\bibitem{Sch_fer_2021}
A.~Schäfer, P.~Chaggar, T.~B. Thompson, A.~Goriely, and E.~Kuhl, ``Predicting brain atrophy from tau pathology: a summary of clinical findings and their translation into personalized models'', {\em Brain Multiphysics}, vol.~2, p.~100039, 2021.

\bibitem{OASIS}
P.~LaMontagne, T.~L. Benzinger, J.~Morris, S.~Keefe, R.~Hornbeck, C.~Xiong, E.~Grant, J.~Hassenstab, K.~Moulder, A.~Vlassenko, M.~Raichle, C.~Cruchaga, and D.~Marcus, ``{OASIS-3: Longitudinal Neuroimaging, Clinical, and Cognitive Dataset for Normal Aging and Alzheimer Disease}'', {\em medRxiv}, 2019.

\bibitem{brainparam1}
E.~Griffiths, J.~Hinrichsen, N.~Reiter, and S.~Budday, ``On the importance of using region-dependent material parameters for full-scale human brain simulations'', {\em European Journal of Mechanics - A/Solids}, vol.~99, p.~104910, 2023.

\bibitem{brainparam2}
J.~Hinrichsen, N.~Reiter, L.~Br\"auer, F.~Paulsen, S.~Kaessmair, and S.~Budday, ``Invaders identification of region-specific hyperelastic material parameters for human brain tissue'', {\em Biomechanics and Modeling in Mechanobiology}, vol.~22, pp.~1729--1749, 2023.

\end{thebibliography}

\end{document}